\documentclass[10pt]{amsart}
\usepackage{cite,mathrsfs,amsmath,amssymb,amsthm,bbm,esint,enumerate,stmaryrd}

\usepackage[breaklinks=true]{hyperref}

\newtheorem{lem}{Lemma}[section]
\newtheorem{thrm}[lem]{Theorem}
\newtheorem{prop}[lem]{Proposition}
\newtheorem{cor}[lem]{Corollary}
\theoremstyle{definition}

\theoremstyle{remark}

\newcommand{\eq}[2]{\begin{equation}\label{#1}#2\end{equation}}

\newcommand{\R}{\mathbb{R}}
\newcommand{\C}{\mathbb{C}}
\newcommand{\Z}{\mathbb{Z}}
\newcommand{\N}{\mathbb{N}}

\newcommand{\Test}{C_c^\infty}

\renewcommand{\Re}{\operatorname{Re}}
\renewcommand{\Im}{\operatorname{Im}}
\renewcommand{\epsilon}{\varepsilon}
\newcommand{\mc}{\mathcal}

\newcommand{\mr}{\mathrm}
\newcommand{\mf}{\mathfrak}

\DeclareMathOperator{\bbE}{\mathbb E}
\DeclareMathOperator{\bbP}{\mathbb P}
\newcommand{\<}{\langle}
\renewcommand{\>}{\rangle}
\newcommand{\p}{\partial}

\newcommand{\bbo}{\mathbbm 1}

\newcommand{\bpf}{\begin{proof}}
\newcommand{\epf}{\end{proof}}
\newcommand{\qtq}[1]{\quad\text{#1}\quad}

\newcommand{\qt}[1]{\quad\text{#1}}

\newcommand{\LHS}[1]{\text{LHS}\eqref{#1}}
\newcommand{\RHS}[1]{\text{RHS}\eqref{#1}}

\newcommand{\cN}{\mc N}
\newcommand{\I}{\mf I}

\newcommand{\loc}{\mr{loc}}
\newcommand{\cI}{\mc I}
\newcommand{\cT}{\mc T}

\newcommand{\opT}{\vert\kern-0.25ex\vert\kern-0.25ex\vert S\vert\kern-0.25ex\vert\kern-0.25ex\vert_T}

\title[Homogenization for the NLS with sprinkled nonlinearity]{Homogenization for the nonlinear Schr\"odinger equation with sprinkled nonlinearity}
\author[B.~Harrop-Griffiths]{Benjamin Harrop-Griffiths}
\address{Department of Mathematics \& Statistics, Georgetown University\\Washington, DC 20057, USA}
\email{benjamin.harropgriffiths@georgetown.edu}
\thanks{}

\author[M.~Ntekoume]{Maria Ntekoume}
\address{Department of Mathematics \& Statistics, Concordia University\\ Montr\'eal, QC H3G 1M8, Canada}
\email{maria.ntekoume@concordia.ca}
\thanks{}

\numberwithin{equation}{section}
\allowdisplaybreaks

\begin{document}

\begin{abstract} We first prove homogenization for the nonlinear Schr\"odinger equation with sprinkled nonlinearity introduced in \cite{2024arXiv240501246H}. We then investigate how solutions fluctuate about the homogenized solution.
\end{abstract}

\maketitle

\section{Introduction}

The nonlinear Schr\"odinger equation with sprinkled nonlinearity,
\eq{sprinkled}{\tag{NLS$_\mu$}
i\p_t\psi = -\p_x^2\psi + 2|\psi|^2\psi\,d\mu,
}
was introduced in \cite{2024arXiv240501246H} as a model for the propagation of waves exhibiting strong self-interaction at defects in the surrounding medium.
These defects are distributed according to a homogeneous Poisson process \(\mu\), which we take to have unit intensity.  Here, and throughout, we use \(d\mu\) to denote the measure \(\mu\) interpreted as a (Schwartz) distribution.  This notation is inspired by the fact that \eqref{sprinkled} describes the flow of the Hamiltonian
\[
\int_\R \tfrac12|\p_x\psi(x)|^2\,dx + \int_\R \tfrac12|\psi(x)|^4\,d\mu(x).
\]

The model \eqref{sprinkled} is motivated by physical scenarios in which nonlinear effects are caused by the interaction of waves with the surrounding environment.
This behavior is illustrated by the optical Kerr effect, in which an intense pulse of light causes a nonlinear change to the refractive index of the material through which it travels (see, e.g., \cite{MR2475397}).

The measure \(\mu\) in \eqref{sprinkled} describes the spatial distribution of the nonlinear effects.
In the case of a homogeneous background, we may take \(\mu\) to be Lebesgue measure and obtain the usual cubic nonlinear Schr\"odinger equation,
\eq{NLS}{\tag{NLS}
i\p_t\psi = - \p_x^2\psi + 2|\psi|^2\psi.
}
At the other extreme, where the nonlinear effects are concentrated at a single point, we may take \(\mu\) to be a Dirac measure. This yields the nonlinear Schr\"odinger equation with concentrated nonlinearity,
\eq{NLS-delta}{\tag{NLS$_\delta$}
i\p_t\psi = - \p_x^2\psi + 2\delta(x)|\psi|^2\psi.
}
The model \eqref{NLS-delta} arises in several physical settings, e.g., \cite{PhysRevB.54.1537,PhysRevB.47.10402,PhysRevB.56.15090,LIDORIKIS1998346,PhysRevE.84.056609,10.1119/1.1417529,F_Kh_Abdullaev_2004}. It has also received considerable attention from the mathematics community; see, e.g., \cite{MR4188177,MR3275343,MR2318828,MR4548487,harropgriffiths2025scatteringnonlinearschrodingerequation,MR1708796,MR1814425,MR3987813,MR4191503,MR4018570,MR4186675}, as well as the review article \cite{MR4653819}.

The sprinkled nonlinearity that we focus on in this paper is a toy model for the scenario in which the nonlinear effects are concentrated at points that have been randomly sprinkled throughout the environment.  The Poisson process \(\mu\) is a natural, physically-motivated model for such a ``random sprinkling''; see \cite[Section 1.3]{2024arXiv240501246H}.
When the typical distance between points at which the nonlinear effects are concentrated is much smaller than a typical wavelength, we might hope to treat the surroundings as if they were homogeneous.
Our goal in this article is to rigorously describe this homogenization procedure.

As discussed in \cite{2024arXiv240501246H}, the choice of a Poisson process in \eqref{sprinkled} is not out of mathematical convenience. Indeed, one may readily envision homogenization occurring when the nonlinearity is concentrated in a periodic arrangement. Taking \(\mu\) to be a Poisson process yields an antithetical scenario: we have regions of arbitrarily high concentration, a potential obstruction to homogenization. Nevertheless, as the nonlinear effects are uniformly distributed on average, we will still be able to prove that homogenization occurs.

The cubic nonlinearity appearing in \eqref{sprinkled} is chosen both for its physical relevance and its convenience.  While we expect that our results can be extended to other sufficiently smooth defocusing nonlinearities with polynomial growth, we have opted for the simplest case of the cubic nonlinearity to avoid additional technical complications.

Given a fixed probability space \((\Omega,\mc F,\bbP)\), we recall (see, e.g., \cite{MR854102}) that a random measure is a measurable function from \(\Omega\) to the space of locally finite Borel measures on \(\R\) endowed with the vague topology.
If \(\mu\) is a random measure and \(E\subseteq \R\) is a bounded Borel set then \(\mu(E)\) is a non-negative random variable.
We say that \(\mu\) has independent increments if \(\mu(E_1),\dots,\mu(E_n)\) are independent random variables for any collection of disjoint bounded Borel sets \(E_1,\dots,E_n\subseteq \R\).
We say that \(\mu\) is stationary if \(\mu(E)\) is identically distributed to \(\mu(E+h)\) for all \(h\in \R\), where \(E+h\) denotes the translation of \(E\) by \(h\).
Finally, we say that \(\mu\) has unit intensity if \(\bbE \mu(E) = m(E)\) for all bounded Borel sets, where \(m\) denotes Lebesgue measure.

The distribution of a random measure is uniquely determined by its Laplace functional, the map
\[
f\mapsto \bbE \exp\left(-\int_\R f\,d\mu\right),
\]
on the set of non-negative \(f\in C_c(\R)
\). It follows from a result of Kingman \cite{MR210185} (also see \cite{MR854102}) that any stationary random measure with independent increments and unit intensity has Laplace functional
\eq{new Phi defn}{
\bbE\exp\left(-\int_\R f\,d\mu\right) = \exp\left(- \int_\R \Phi(f) + f\,dx\right),
}
where 
\eq{new Phi expression}{
\Phi(z) = \int_0^\infty 1 - sz - e^{-sz}\,d\Lambda(s),
}
for a Borel measure \(\Lambda\) on \((0,\infty)\) satisfying \(\int_0^\infty s\,d\Lambda(s)<\infty\).  We call \(\Lambda\) the L\' evy measure of \(\mu\).

A homogeneous Poisson process of unit intensity is a stationary random measure \(\mu\) with independent increments so that for any bounded Borel set \(E\subseteq\R\) the random variable \(\mu(E)\) is Poisson-distributed with mean \(m(E)\).
This corresponds to taking \(\Lambda\) to be a Dirac mass at \(s=1\) in \eqref{new Phi expression}, so that
\eq{Poisson Phi}{
\Phi(z) = 1 - z - e^{-z}.
}

For the remainder of this article, we assume that \(\mu\) is a stationary random measure with independent increments and unit intensity so that the associated L\' evy measure \(\Lambda\) satisfies
\eq{Levy exponential bound}{
\int_0^\infty s\,e^{as}\,d\Lambda(s)<\infty,
}
for some \(a>0\). This ensures that \(\Phi(z)\) is analytic for \(\Re z>-a\) and satisfies \(\Phi(0) = \Phi'(0) = 0\). By rescaling, we may also assume that \(\Phi''(0) = -1\). We consider the homogeneous Poisson process to be the canonical example of such a random measure; other examples include compound Poisson processes and Gamma processes.

To investigate the problem of homogenization, we take a sequence of positive numbers \(0<\epsilon_n\to0\) and consider a corresponding sequence of solutions \(\psi_n\colon \R\times \R\to \C\) of the equation
\eq{NLSn}{\tag{NLS$_n$}
i\p_t\psi_n = -\p_x^2\psi_n + 2|\psi_n|^2\psi_n\,d\mu_n,
}
where \(\mu_n\) is the unique random measure with Laplace functional
\begin{equation}\phantomsection\label{Laplace functional}
\bbE\exp\biggl(-\int_\R f\,d\mu_n\biggr) = \exp\biggl(- \tfrac1{\epsilon_n}\int_\R\Phi(\epsilon_nf) + \epsilon_n f\,dx\biggr),
\end{equation}
for non-negative \(f\in C_c(\R)\).
The equation \eqref{NLSn} describes the large scale behavior of \eqref{sprinkled}: the rescaled random measure \(E\mapsto\frac1{\epsilon_n}\mu_n(\epsilon_nE)\) is identically distributed to \(\mu\). Consequently, \(\psi_n(t,x)\) solves \eqref{NLSn} if and only if
\[
\psi(t,x) = \epsilon_n\psi_n(\epsilon_n^2t,\epsilon_n x)
\]
solves \eqref{sprinkled}.
It follows\footnote{While \cite{2024arXiv240501246H} only considers the case that \(\mu\) is a homogeneous Poisson process of unit intensity, the only place that any specific properties of the random measure are used is \cite[Lemma 4.1]{2024arXiv240501246H}. One may readily verify that the proof of this lemma applies without modification to the more general random measures \(\mu\) that we consider.  In fact, all one requires is that \(\sup_{\ell\in \Z}\bbE\mu\bigl([\ell-\frac12,\ell+\frac12)\bigr)^p<\infty\) for all sufficiently large \(p\geq 1\); see Lemma~\ref{l:this one} below for a slight generalization.} from \cite[Theorem 1.2]{2024arXiv240501246H} that the equation \eqref{NLSn} is almost surely globally well-posed in \(H^1\).

Our first result proves that the solutions \(\psi_n\) of \eqref{NLSn} homogenize to the solution \(\psi\) of \eqref{NLS}:

\begin{thrm}[Homogenization]\label{t:homogenization}
Let \(0<\epsilon_n\to0\) and \(\psi_0\in H^1\). If \(\psi_n\in C(\R;H^1)\) is the almost surely defined solution of \eqref{NLSn} with initial data \(\psi_n(0) = \psi_0\) and \(\psi\in C(\R;H^1)\) is the solution of \eqref{NLS} with initial data \(\psi(0) = \psi_0\) then, for all \(T>0\) and \(1\leq p<\infty\), we have \(\psi_n\to \psi\) in \(L^p\bigl(d\bbP;C\bigl([-T,T];H^1\bigr)\bigr)\), i.e.,
\eq{cvgce}{
\lim_{n\to\infty}\bbE\sup_{|t|\leq T}\|\psi_n(t) - \psi(t)\|_{H^1}^p = 0.
}
\end{thrm}

Theorem~\ref{t:homogenization} tells us that, for a fixed choice of initial data and compact time interval, the solution \(\psi_n\) of \eqref{NLSn} is well-approximated on average by the solution \(\psi\) of \eqref{NLS}, provided \(n\gg1\) is sufficiently large. A natural next question is whether we can describe how the solutions \(\psi_n\) fluctuate about the homogenized solution \(\psi\).

Recalling that we denote Lebesgue measure by \(m\), we may formally write
\[
\mu_n = m + \sqrt{\epsilon_n}\nu_n,\qtq{with}\nu_n = \tfrac1{\sqrt{\epsilon_n}}(\mu_n - m).
\]
Considering \(\nu_n\) to be a (Schwartz-)distribution-valued random variable, we may use \eqref{Laplace functional} to see that \(\nu_n\) converges in distribution to white noise, the real-Schwartz-distribution-valued random variable \(\xi\) for which
\eq{WN defn}{
\bbE \exp\bigl(i\<f,\xi\>\bigr) = \exp\bigl(-\tfrac12\|f\|_{L^2}^2\bigr),
}
for all real-valued \(f\in \Test(\R)\).
Formally expanding \eqref{NLSn} in powers of \(\sqrt{\epsilon_n}\) now shows us that
\eq{varphin}{
\varphi_n = \tfrac1{\sqrt{\epsilon_n}}(\psi_n - \psi)
}
should converge (in some sense) to a solution of the equation
\eq{lin-NLS}{
i\p_t\varphi = - \p_x^2\varphi 
+ 4 |\psi|^2\varphi + 2\psi^2\bar\varphi + 2|\psi|^2\psi\,\xi,
}
with initial data \(\varphi(0) = 0\).

We emphasize here that \(\xi\) is \emph{spatial} white noise.  Consequently, the analysis of \eqref{lin-NLS} reduces to understanding the linearization of \eqref{NLS} about the solution \(\psi\).  This is significantly more straightforward than the case of temporal noise considered in, e.g., \cite{MR3236753}.  In particular, for all \(\psi\in C(\R;H^1)\) and \(s>\frac12\) there almost surely exists a unique solution \(\varphi\in C(\R;H^{-s})\) of \eqref{lin-NLS} satisfying \(\varphi(0) = 0\).  Further, for all \(t\in \R\) the random function \(\varphi(t)\in H^{-s}\) is a mean zero Gaussian process with covariance and pseudo-covariance determined by \(t\) and \(\psi\); see Proposition~\ref{p:WP of lin-NLS}.

Our second main result proves the convergence of \(\varphi_n\) to \(\varphi\):

\begin{thrm}[Fluctuations]\label{t:cvgce in distrib}
Let \(s>\frac12\), \(0<\epsilon_n\to0\), and \(\psi_0\in H^1\). Let \(\psi_n\in C(\R;H^1)\) be the almost surely defined solution of \eqref{NLSn} with initial data \(\psi_n(0) = \psi_0\) and \(\psi\in C(\R;H^1)\) be the solution of \eqref{NLS} with initial data \(\psi(0) = \psi_0\). Let \(\varphi_n\) be defined as in \eqref{varphin} and \(\varphi\in C(\R;H^{-s})\) be the almost surely defined solution of \eqref{lin-NLS} with initial data \(\varphi(0) = 0\). Then, for any \(T>0\) the sequence \(\varphi_n\) converges to \(\varphi\) in distribution as \(C\bigl([-T,T];H^{-s}\bigr)\)-valued random variables.
\end{thrm}

Given that \(\varphi\) is a Gaussian process, Theorem~\ref{t:cvgce in distrib} demonstrates that the solutions \(\psi_n\) of \eqref{NLSn} exhibit Gaussian fluctuations about their mean.  As the covariance and pseudo-covariance of \(\varphi\) can be explicitly written in terms of \(\psi\), one might hope to be able to give a detailed description of the statistics of these fluctuations.  Unfortunately, the expressions for the covariance and pseudo-covariance of \(\varphi\) derived in Proposition \ref{p:WP of lin-NLS} are rather opaque. We leave it as an interesting open problem to describe these statistics in more detail.

Homogenization problems in PDE abound, especially in the context of elliptic equations.
Early works, such as \cite{MR348255}, focused on homogenization for deterministic linear elliptic operators with periodic coefficients.
The study of stochastic homogenization, where the operator is taken to have random coefficients, was initiated shortly thereafter in \cite{MR712714,MR542557}. Nonlinear problems were subsequently considered in \cite{MR850613,MR870884}.
We refer the reader to the recent article \cite{armstrong2022elliptichomogenizationqualitativequantitative} for an introduction to elliptic homogenization and a discussion of known results.

Homogenization problems for (time-dependent) Schr\"odinger-type equations have largely focused on linear models. As in the elliptic case, authors have considered both periodic homogenization \cite{MR2839402,MR2166838,MR1421214,MR4683124,MR4698781,MR2530894,MR3646722,MR3563045,MR4484681,MR4039349} and stochastic homogenization \cite{MR3159467,MR3187780}.

By contrast, homogenization for nonlinear Schr\"odinger equations has received relatively little attention.
In \cite{MR2216722}, the author considered a nonlinear Schr\"odinger equation with a periodic potential and proved homogenization for well-prepared initial data, providing a nonlinear analog of the effective mass theorems in \cite{MR2166838,MR1421214}.

To the best of our knowledge, homogenization questions of the type considered in this article, where the inhomogeneity arises from the nonlinear interactions, have only been considered in previous work of the second author \cite{MR4162949}. There, subtle perturbation techniques in critical spaces were used to derive a sufficient condition for the coupling function \(g(x)\) to ensure homogenization of the equation
\[
i\p_t\psi_n = - \Delta\psi_n + g(nx)|\psi_n|^2\psi_n
\]
in \(L^2(\R^2)\).
While the sufficient condition obtained in \cite{MR4162949} was general enough to encompass several physically motivated examples, both deterministic and random, an important limitation was that it was only applicable to bounded coupling functions, usually of periodic nature. Achieving homogenization in the unbounded setting, when defects are allowed to concentrate, is one of the key contributions of this paper.

Given that Theorems~\ref{t:homogenization} and~\ref{t:cvgce in distrib} apply to a reasonably large class of random measures, the homogenized equation \eqref{NLS} and the equation for the fluctuations \eqref{lin-NLS} are, in a sense, universal. 
With this perspective, one may view Theorem~\ref{t:cvgce in distrib} as a central-limit-type theorem. Similar questions about the structure of random perturbations in homogenization problems have been studied in both the elliptic and parabolic contexts (see, e.g., \cite{MR1728027,MR2925905,MR2435670,MR3554429}) but never, as far as we are aware, in the Schr\"odinger case.

An important question is whether we can expand the class of random measures to which Theorems~\ref{t:homogenization} and~\ref{t:cvgce in distrib} apply even further.
An inspection of our proof demonstrates that the main bottleneck arises in Proposition~\ref{p:Haar}. Here, we use that the measure \(\mu\) is stationary and has independent increments to derive estimates for a collection of orthonormal random variables \(X_{N,k}\).  It seems reasonable to conjecture that one could adapt the proof to use an appropriate notion of almost orthogonality, which would make it possible to consider more general measures \(\mu\).  We leave this as an open problem to be considered in future work.

The (deterministic) model \eqref{NLS-delta} arises from physical settings in which the nonlinear effects are concentrated near \(x=0\) at much smaller scales than typical wavelengths.  To this end, it was proved in \cite{MR3275343,MR2318828,harropgriffiths2025scatteringnonlinearschrodingerequation} that the dynamics of \eqref{NLS-delta} can be obtained as the \(h\to0\) limit of a mollified model, where the distribution \(\delta(x)\) is replaced by the function \(\zeta^h(x) = \frac1h\zeta\bigl(\frac xh\bigr)\) for a non-negative, even function \(\zeta\in \Test(-1,1)\) satisfying \(\int_\R\zeta\,dx = 1\). A corresponding result for the random measure \(\mu\) was proved in \cite{2024arXiv240501246H}. Indeed, it follows from \cite[Theorem 1.3]{2024arXiv240501246H} that, if \(\psi_n^h\) is the (almost surely defined) solution of \eqref{NLSn} with \(\mu_n\) replaced by the absolutely continuous measure \(\mu_n^h\) with Radon--Nikodym derivative
\eq{RND}{
\frac{d\mu_n^h}{dx} = \int_\R \zeta^h(x-y)\,d\mu_n(y),
}
we have \(\psi_n^h\to\psi_n\) in \(L^p\bigl(d\bbP;C\bigl([-T,T];H^1\bigr)\bigr)\) as \(h\to 0\).  In Appendix~\ref{app:mollification} we prove that it is possible to obtain analogs of both Theorems~\ref{t:homogenization} and~\ref{t:cvgce in distrib} for the mollified measures \(\mu_n^h\) that hold uniformly for \(0<h\leq 1\).

Let us now outline the proof of our main results. Given a non-negative sequence \(Z\colon \Z\to[0,\infty)\), we denote
\eq{LZ2}{
\|f\|_{L_Z^2}^2 = \int_\R |f(x)|^2\,\omega(x;Z)\,dx,
}
where the piecewise linear weight
\eq{def: omega}{
\omega(x;Z) = \sum_{k\in \Z}\Bigl[(1+k-x)\,|\cN(k;Z)|^2 + (x-k)\,|\cN({k+1};Z)|^2\Bigr]\bbo_{[k,k+1)}(x),
}
and the slowly varying envelope \(\cN(k;Z)>0\) for the sequence \(Z\) is given by
\eq{defn: cnk}{
|\cN(k;Z)|^2 = 4 + \sup_{\ell\in \Z}\Bigl[|Z(\ell)|^2 - |k - \ell|\Bigr].
}

We will use this construction with two particular choices of the sequence \(Z\). Precisely, for each sample of our random measure \(\mu_n\) we take \(X_n^1\subseteq H^1\) to consist of all \(f\in H^1\) with finite norm
\eq{Xn1}{
\|f\|_{X_n^1}^2 = \|f\|_{H^1}^2 + \|f\|_{L_{Z_n}^2}^2,\qtq{where}
Z_n(k) = \mu_n\bigl([k-\tfrac12,k+\tfrac12)\bigr),
}
and for \(\frac12<s\leq 1\) we take \(Y_{n,s}^1\subseteq X_n^1\) to consist of all \(f\in X_n^1\) with finite norm
\eq{Yns1}{
\|f\|_{Y_{n,s}^1}^2 = \|f\|_{X_n^1}^2 + \|f\|_{L_{Z_{n,s}}^2}^2,\qtq{where}Z_{n,s}(k) = \bigl\|\rho_k\tfrac{d\mu_n - 1}{\sqrt{\epsilon_n}}\bigr\|_{H^{-s}},
}
and \(1 = \sum_{k\in \Z}\rho_k\) is a partition of unity so that $\rho_k(x)=\rho(x-k)$ for a nonnegative function $\rho\in C_c^\infty (-1,1)$.
Recalling that we distinguish between the interpretation of \(\mu_n\) as a measure and as a distribution by writing \(d\mu_n\) in the latter case, for \(f\in \Test(\R)\) we have
\[
\<f,d\mu_n - 1\> = \int_\R f\,d\mu_n - \int_\R f\,dx.
\]

For a fixed sample of the random measure \(\mu_n\), the results of \cite{2024arXiv240501246H} prove that the equation \eqref{NLSn} is globally well-posed in \(X_n^1\); see Lemma~\ref{l:ap psin} below. Crucially, we will prove in Corollary~\ref{c:L2mu exp 1} that if \(\psi_0\in H^1\) then \(\psi_0\in \bigcap_{n\geq 1}X_n^1\) almost surely and hence the sequence of solutions \(\psi_n\) is almost surely well-defined.

To prove Theorem~\ref{t:homogenization}, we start by expressing the difference
\begin{align}\label{intro difference}
\psi_n(t) - \psi(t) &= -2i\int_0^t e^{i(t-\tau)\p_x^2}\Bigl[\bigl(|\psi_n(\tau)|^2\psi_n(\tau) - |\psi(\tau)|^2\psi(\tau)\bigr)\,d\mu_n\Bigr]\,d\tau\\
&\quad - 2i\int_0^te^{i(t - \tau)\p_x^2}\Bigl[|\psi(\tau)|^2\psi(\tau) \,\bigl(d\mu_n-1\bigr)\Bigr]\,d\tau.\notag
\end{align}
We then consider a similar circle of estimates to those used in the proof of well-posedness in \cite{2024arXiv240501246H}; see also \cite{MR2318828,MR3275343}. Precisely, we first bound the difference \(\psi_n(t) - \psi(t)\) in \(L^\infty\), then in \(H^{-1}\), and then use conservation of mass and energy to bound the difference of the \(H^1\)-norms of \(\psi_n(t)\) and \(\psi(t)\).
Finally, we upgrade to convergence in \(H^1\) using the Radon--Riesz Theorem.

While the first term on \RHS{intro difference} may be controlled using similar estimates to \cite{2024arXiv240501246H}, the second term requires us to estimate \(\psi(t)\) in the new space \(Y_{n,s}^1\).
To prove this estimate, we must establish two key facts.

The first is that the argument used to obtain a priori estimates in \cite{2024arXiv240501246H} can be extended to estimate \(\psi\) in the space \(L^p\bigl(d\bbP;C\bigl([-T,T];L_Z^2\bigr)\bigr)\) for any random sequence \(Z\) satisfying the uniform bound
\[
\sup_{k\in \Z}\bbE |Z(k)|^{6p}\lesssim_p1.
\]
This fact is behind the a priori estimates obtained in Corollary~\ref{c:L2mu exp 2} and Lemma~\ref{l:ap psi}.

The second fact we require is that the moments of the sequence \(Z_{n,s}(k)\) can be bounded uniformly in both \(n\) and \(k\). This is proved in Proposition~\ref{p:Haar} by considering a sequence of centered random variables,
\[
X_{N,k} = \int_\R e_{N,k}\,d\mu_n - \int_\R e_{N,k}\,dx,
\]
where \(e_{N,k}\) is a Haar basis of \(L^2\). By computing the joint cumulants of products of these random variables, we can apply the moment-cumulant formula to estimate the moments.

Turning to the proof of Theorem~\ref{t:cvgce in distrib}, we first note that Proposition~\ref{p:WP of lin-NLS} follows from the well-posedness of the linearization of \eqref{NLS} about the solution \(\psi\).

To prove that \(\varphi_n\) converges to \(\varphi\) in distribution, we must prove that for all bounded continuous functions \(F\colon C([-T,T];H^{-s})\to \C\) we have
\[
\bbE F(\varphi_n)\to \bbE F(\varphi),
\]
or equivalently that the corresponding probability distributions,
\[
(\varphi_n)_*\bbP(E) = \bbP\bigl[\varphi_n\in E\bigr],
\]
converge weakly to \(\varphi_*\bbP\) in the space of Borel probability measures on\linebreak \(C([-T,T];H^{-s})\). To achieve this, we will first prove (see Proposition~\ref{p:compactness}) that our sequence of probability measures is tight. Prokhorov's Theorem then ensures that the sequence is weakly (sequentially) precompact. To complete the proof we will show that for every bounded linear functional \(\cT\colon C([-T,T];H^{-s})\to \C\) the characteristic functional satisfies
\eq{fluctuation goal}{
\lim_{n\to\infty}\bbE e^{i\Re \cT(\varphi_n)} = \bbE e^{i\Re \cT(\varphi)}.
}
This proves that the only weak limit point of the sequence of measures \((\varphi_n)_*\bbP\) is \(\varphi_*\bbP\), thereby completing the proof.

The main technical ingredient in this argument is a collection of a priori estimates for \(\varphi_n\) under the hypothesis that both \(\psi_n\) and \(\psi\) are bounded in \(C\bigl([-T,T];Y_{n,s}^1\bigr)\) by a given constant. These follow from refining the \(H^{-1}\) estimates for \eqref{intro difference} used in the proof of Theorem~\ref{t:homogenization}.
Crucially, we now require estimates for \(\psi_n\) in\linebreak \(L^p\bigl(d\bbP;C\bigl([-T,T];Y_{n,s}^1\bigr)\bigr)\), not just \(\psi\). Fortunately, our proof of the corresponding estimate for \(\psi\) is sufficiently robust that it can be applied to \(\psi_n\) with essentially no modifications.

The remainder of this article is structured as follows: in Section~\ref{s:prelim} we introduce some notation and prove several preliminary estimates. This includes the key estimate \eqref{Haar moments} that yields uniform bounds for the moments of \(Z_{n,s}(k)\). Section~\ref{s:homogenization} is dedicated to the proof of Theorem~\ref{t:homogenization}, with the main estimate for the difference \(\psi_n - \psi\) appearing in Proposition~\ref{p:the beast}. Finally, in Section~\ref{s:fluctuations} we first prove Proposition~\ref{p:WP of lin-NLS} and then Theorem~\ref{t:cvgce in distrib}. We conclude with Appendix~\ref{a:linear}, in which we discuss the construction and properties of the linear propagator for the linearization of \eqref{NLS}, and Appendix~\ref{app:mollification}, in which we discuss the adaptation of Theorems~\ref{t:homogenization} and~\ref{t:cvgce in distrib} to the mollified measures \(\mu_n^h\).

\subsection*{Acknowledgments} BHG was supported by NSF grant DMS-2406816. MN was supported by NSERC grant RGPIN-2024-05082. The authors are grateful to Bjoern Bringmann for some informative comments regarding an early version of the proof. The authors are also indebted to the anonymous referees for their numerous helpful comments and suggestions.

\section{Preliminaries}\label{s:prelim}

We start by introducing some notation. Given a Banach space \(X\) of functions on \(\R\), a time \(T>0\), and \(f\in L^p\bigl((-T,T);X\bigr)\) for \(1\leq p< \infty\), we denote the norm
\[
\|f\|_{L_T^p X} = \biggl[\int_{-T}^T \|f(t)\|_X^p\,dt\biggr]^{\frac1p}.
\]
We extend this notation to continuous functions \(f\in C\bigl([-T,T];X\bigr)\) by writing
\[
\|f\|_{C_TX} = \sup_{|t|\leq T}\|f(t)\|_X.
\]

We take \(\chi\in \Test(\R)\) to be a fixed even function satisfying \(0\leq \chi\leq 1\) that is identically \(1\) on \([-1,1]\) and supported on \((-2,2)\). Given \(R>0\), we rescale to obtain \(\chi_R(x) = \chi\bigl(\tfrac xR\bigr)\).

We normalize the Fourier transform by setting
\[
\hat f(\xi) = \tfrac1{\sqrt{2\pi}}\int_\R f(x)e^{-ix\xi}\,dx,
\]
so that for \(s\in \R\) we have
\[
\|f\|_{H^s}^2 = \int_\R \<\xi\>^{2s}|\hat f(\xi)|^2\,d\xi,
\]
where \(\<\xi\>^2 = 1 + |\xi|^2\). Given dyadic \(N\geq 1\), we denote the Littlewood--Paley projection
\[
P_{\leq N}f = \tfrac1{\sqrt{2\pi}}\int_\R \chi_N(\xi)\hat f(\xi)e^{ix\xi}\,d\xi,
\]
and set \(P_{>N} = 1 - P_{\leq N}\).

We recall that for all \(s>\frac12\) we have \(H^s(\R)\subseteq C_0(\R)\subseteq L^\infty(\R)\), with the Gagliardo--Nirenberg inequality
\eq{GN}{
\|f\|_{L^\infty}\lesssim \|f\|_{L^2}^{1-\frac1{2s}}\|f\|_{H^s}^{\frac1{2s}}.
}

If \(k\geq 0\) is an integer, we may apply the Leibniz rule to bound
\[
\|fg\|_{H^k}\lesssim \|f\|_{W^{k,\infty}}\|g\|_{H^k}.
\]
By interpolation and duality, we may then show that for all \(|s|\leq k\) we have the product estimate
\eq{product estimate inf}{
\|fg\|_{H^s}\lesssim \|f\|_{W^{k,\infty}}\|g\|_{H^s}.
}

A second product estimate that we will use frequently is:
\begin{lem}
For all \(|s|\leq 1\) we have
\eq{product estimate}{
\|fg\|_{H^s}\lesssim \|f\|_{H^1}\|g\|_{H^s}.
}
\end{lem}
\bpf
For \(s=0\), we apply \eqref{GN} to bound
\[
\|fg\|_{L^2}\leq \|f\|_{L^\infty}\|g\|_{L^2}\lesssim \|f\|_{H^1}\|g\|_{L^2},
\]
and for \(s=1\) we combine this with the estimate
\[
\|(fg)_x\|_{L^2}\lesssim \|f_x\|_{L^2}\|g\|_{L^\infty} + \|f\|_{L^\infty}\|g_x\|_{L^2}\lesssim \|f\|_{H^1}\|g\|_{H^1}.
\]

For \(0<s<1\), we choose \(2<p,q<\infty\) so that \(\frac12 = \frac1p+\frac1q\) and \(s-\frac12\leq \frac1p\leq s\). We then apply the fractional Leibniz rule to bound
\[
\|fg\|_{H^s}\lesssim \|f\|_{L^\infty}\|g\|_{H^s} + \|\<\p_x\>^sf\|_{L^p}\|g\|_{L^q}.
\]
The first summand is readily bounded using \eqref{GN}. For the second, we apply the Hardy--Littlewood--Sobolev inequality to obtain
\[
\|\<\p_x\>^s f\|_{L^p}\lesssim \|f\|_{H^{s+\frac1q}}\lesssim \|f\|_{H^1}\qtq{and}\|g\|_{L^q}\lesssim \|g\|_{H^{\frac1p}}\lesssim \|g\|_{H^s}.
\]

Finally, for the case \(-1\leq s<0\) we argue by duality. Given a test function \(h\in H^{-s}\), we apply the result for positive \(s\) to estimate
\[
|\<fg,h\>|\lesssim \|g\|_{H^s}\|\bar fh\|_{H^{-s}}\lesssim \|g\|_{H^s}\|f\|_{H^1}\|h\|_{H^{-s}},
\]
which completes the proof of \eqref{product estimate}.
\epf

Our next lemma characterizes precompactness in the space \(C_TH^s\):
\begin{lem}\label{l:compactness helper}
Let \(T>0\) and \(s\in \R\). A bounded set \(\mc F\subseteq C\bigl([-T,T];H^s\bigr)\) is precompact if and only if it is uniformly equicontinuous in time,
\eq{time equity}{
\lim_{h\to 0}\sup_{\substack{|t-\tau|\leq h\\t,\tau\in [-T,T]}}\sup_{f\in \mc F}\|f(t) - f(\tau)\|_{H^s} = 0,
}
and both \(H^s\)-tight and \(H^s\)-equicontinuous in space, uniformly in time,
\begin{align}\label{compactness helper}
\lim_{R\to \infty}\sup_{f\in \mc F}\|(1 - \chi_R)f\|_{C_TH^s} = 0 = \lim_{N\to\infty}\sup_{f\in \mc F}\|P_{>N}f\|_{C_TH^s}.
\end{align}
\end{lem}
\bpf
Let us start by recalling from \cite[Theorem 3]{MR0801333} that a bounded set \(Q\subseteq L^2\) is precompact if and only if
\eq{compactness helper half}{
\lim_{R\to\infty}\sup_{q\in Q}\|(1 - \chi_R)q\|_{L^2} = 0 = \lim_{N\to\infty}\sup_{q\in Q}\|P_{>N}q\|_{L^2}.
}

To extend this to \(H^s\), we first use Plancherel's Theorem to prove that for all \(R\geq 1\) we have
\begin{align}
&\Bigl\|\bigl[\<\p_x\>^s,\chi_R\bigr] q\Bigr\|_{L^2}\notag\\
&\qquad\simeq \left\|\int_\R \hat\chi_R(\xi-\eta)\bigl[\<\xi\>^s - \<\eta\>^s\bigr]\hat q(\eta)\,d\eta\right\|_{L_\xi^2}\notag\\
&\qquad\simeq \left\|\int_\R\hat \chi_R(\xi - \eta)\left[\int_0^1 s\<\eta + h(\xi-\eta)\>^{s-2}(\eta + h(\xi-\eta))(\xi-\eta)\,dh\right]\hat q(\eta)\,d\eta\right\|_{L_\xi^2}\notag\\
&\qquad\lesssim  \|\<\xi\>^{|s-1|}\xi\,\hat \chi_R(\xi)\|_{L_\xi^1}\|\<\eta\>^{s-1}\hat q(\eta)\|_{L^2_\eta}\lesssim \tfrac1R\|q\|_{H^{s-1}},\label{something}
\end{align}
where we have used that for all \(\eta,\zeta\in \R\) we have
\[
\<\eta + \zeta\>^{s-2}|\eta + \zeta|\leq \<\eta + \zeta\>^{s-1}\lesssim \<\zeta\>^{|s-1|}\<\eta\>^{s-1}.
\]
Indeed, when \(|\zeta|\ll|\eta|\) we have
\[
\<\eta + \zeta\>^{s-1}\simeq\<\eta\>^{s-1}\lesssim \<\zeta\>^{|s-1|}\<\eta\>^{s-1},
\]
whereas when \(|\eta|\lesssim |\zeta|\) we have
\[
\<\eta + \zeta\>^{s-1}\lesssim \<\zeta\>^{\max\{s-1,0\}}\lesssim \<\zeta\>^{|s-1|}\<\eta\>^{s-1}.
\]

Applying \eqref{compactness helper half} to \(\<\p_x\>^sq\) with our commutator estimate, we now see that a bounded set \(Q\subseteq H^s\) is precompact if and only if
\eq{compactness helper half s}{
\lim_{R\to\infty}\sup_{q\in Q}\|(1 - \chi_R)q\|_{H^s} = 0 = \lim_{N\to\infty}\sup_{q\in Q}\|P_{>N}q\|_{H^s}.
}

By the (generalized) Arzel\`a--Ascoli Theorem, a set \(\mc F\subseteq C\bigl([-T,T];H^s\bigr)\) is precompact if and only if it satisfies \eqref{time equity} and for each \(t\in [-T,T]\) the set
\(
\bigl\{f(t):f\in \mc F\bigr\}\subseteq H^s
\)
is precompact. In light of \eqref{compactness helper half s}, it remains to show that if \(\mc F\subseteq C\bigl([-T,T];H^s\bigr)\) is precompact then we have the uniform-in-time estimates \eqref{compactness helper}.

Given \(\eta>0\), we choose an \(\eta\)-covering \(\{f_1,\dots,f_L\}\subseteq \mc F\) of \(\mc F\) and take
\[
Q = \bigcup_{\ell=1}^L\bigl\{f_\ell(t):t\in[-T,T]\bigr\},
\]
which is a finite union of compact subsets of \(H^s\) and therefore compact. Applying \eqref{compactness helper half s} then yields
\[
\lim_{R\to\infty}\max_{\ell=1,\dots,L}\|(1 - \chi_R)f_\ell\|_{C_TH^s} = 0 = \lim_{R\to\infty}\max_{\ell=1,\dots,L}\|P_{>N}f_\ell\|_{C_TH^s}.
\]
Given any \(f\in \mc F\), we may then estimate
\begin{align*}
\|(1 - \chi_R)f\|_{C_TH^s}&\leq \max_{\ell=1,\dots,L}\|(1 - \chi_R)f_\ell\|_{C_TH^s} + \eta,\\
\|P_{>N}f\|_{C_TH^s}&\leq \max_{\ell=1,\dots,L}\|P_{>N}f_\ell\|_{C_TH^s} + \eta,
\end{align*}
where we note that the first of these uses \eqref{product estimate inf}. The proof is now completed by taking the supremum over \(f\in \mc F\) and then sending \(R\to\infty\), respectively \(N\to\infty\).
\epf

In the proof of Theorem~\ref{t:cvgce in distrib}, we will use the following corollary to Lemma~\ref{l:compactness helper}:

\begin{cor}\label{c:FM}
Let \(T>0\) and \(s\in \R\). Given sequences of dyadic integers \(1\leq N_j,R_j\to\infty\) and a constant \(M>0\), let \(\mc F_M\) consist of all \(f\in C_TH^s\) so that
\begin{align}
\|f\|_{C_TH^s} + \|\p_tf\|_{C_TH^{s-2}} &\leq M,\label{FM1-i}\\
\sum_{j=1}^\infty 2^j\|P_{>N_j}f\|_{C_TH^{s}}&\leq M,\label{FM2-i}\\
\sum_{j=1}^\infty 2^j\|(1 - \chi_{R_j})f\|_{C_TH^{s-1}} &\leq M.\label{FM3-i}
\end{align}
Then \(\mc F_M\) is a precompact subset of \(C_TH^{s}\).
\end{cor}
\bpf
We first note that \eqref{FM1-i} ensures that \(\mc F_M\) is bounded in \(C_TH^{s}\) and \eqref{FM2-i} shows that it is equicontinuous in \(H^{s}\), uniformly in time,
\eq{main step for equicont-i}{
\lim_{N\to\infty} \sup_{f\in \mc F_M}\|P_{>N}f\|_{C_TH^{s}} = 0.
}
Given any \(t,\tau\in[-T,T]\), we may apply Bernstein's inequality to bound
\begin{align*}
\|f(t) - f(\tau)\|_{H^{s}} & \leq \bigl\|P_{\leq N}\bigl[f(t) - f(\tau)\bigr]\bigr\|_{H^{s}} + 2\|P_{>N}f\|_{C_TH^{s}}\\
&\lesssim N^2\|\p_tf\|_{C_TH^{s-2}}|t - \tau| + \|P_{>N}f\|_{C_TH^{s}}.
\end{align*}
Taking first the supremum over \(f\in \mc F_M\), then the limit as \(|t - \tau|\to 0\), and finally the limit as \(N\to\infty\) using \eqref{main step for equicont-i}, we see that \(\mc F_M\) is uniformly equicontinuous in time.
Similarly, for all \(R\geq 1\) we may bound
\begin{align*}
\|(1 - \chi_R)f\|_{C_TH^{s}}&\lesssim N\|(1 - \chi_R)f\|_{C_TH^{s-1}} + \|P_{>N} f\|_{C_TH^{s}}.
\end{align*}
Again taking the supremum over \(f\in \mc F_M\), the limit as \(R\to\infty\) using \eqref{FM3-i}, and finally the limit as \(N\to\infty\) using \eqref{main step for equicont-i} shows that \(\mc F_M\) is tight in \(H^{s}\), uniformly in time. We may now apply Lemma~\ref{l:compactness helper} to conclude that \(\mc F_M\) is precompact in \(C_TH^{s}\).
\epf

Let us now turn to the weighted space \(L_Z^2\) defined in \eqref{LZ2}. By construction (see \cite[Lemma 2.1]{2024arXiv240501246H}), we have
\[
\Bigl||\cN(k;Z)|^2 - |\cN(k+1;Z)|^2\Bigr|\leq 1,
\]
and hence the map \(x\mapsto \omega(x;Z)\) is \(1\)-Lipschitz. As a consequence, we have the estimate
\eq{Lipschitz bound}{
\omega(x;Z)\leq |\cN(0;Z)|^2 + |x|,
}
so if \(\cN(0;Z)<\infty\) then \(\cN(k;Z) = \sqrt{\omega(k;Z)}<\infty\) for all \(k\in \Z\).

By definition, we have \(Z(k)\leq \cN(k;Z)\) and we may readily verify that \(\cN(k;Z)\simeq \cN(\ell;Z)\) whenever \(|k - \ell|\leq 3\), where the implicit constants are independent of \(Z\). As a consequence, if \(1 = \sum_{k\in \Z}\rho_k\) is the partition of unity appearing in \eqref{Yns1} and \(\cN(0;Z)<\infty\), then
\eq{sq sum}{
\|f\|_{L_Z^2}^2 \simeq \sum_{k\in \Z} |\cN(k;Z)|^2 \|\rho_k f\|_{L^2}^2.
}
It will also be useful to note that
\eq{sq sum H1}{
\|f\|_{L^2}^2\simeq\sum_{k\in \Z}\|\rho_k f\|_{L^2}^2\qtq{and}\|f\|_{H^1}^2\simeq \sum_{k\in \Z}\|\rho_k f\|_{H^1}^2.
}

We subsequently restrict our attention to the measures \(\mu_n\) with Laplace functionals \eqref{Laplace functional}, and take \(Z_n\), \(Z_{n,s}\) to be the associated random sequences defined in \eqref{Xn1}, \eqref{Yns1}.

We then have the following adaptation of \cite[Lemma 2.3]{2024arXiv240501246H}:
\begin{lem}
If \(Z_n\) is defined as in \eqref{Xn1} and  \(\cN(0;Z_n)<\infty\) then for all \(f_1,f_2\in X_n^1\) and \(g\in C_0(\R)\) we have the estimates
\begin{align}
\int_\R |f_1f_2g|\,d\mu_n & \lesssim \|f_1\|_{X_n^1}\|f_2\|_{X_n^1}\|g\|_{L^\infty},\label{technical 1}\\
\|f_1f_2g\,d\mu_n\|_{H^{-1}} & \lesssim \|f_1\|_{X_n^1}\|f_2\|_{X_n^1}\|g\|_{L^\infty}.\label{technical 2}
\end{align}

Further, if \(\frac12<s\leq 1\) and \(Z_{n,s}\) is defined as in \eqref{Yns1} then whenever \(\cN(0;Z_n)<\infty\), \(\cN(0;Z_{n,s})<\infty\), and \(f_1,f_2,g\in Y_{n,s}^1\) we have the estimate
\begin{align}
\bigl\|f_1f_2g\tfrac{d\mu_n-1}{\sqrt{\epsilon_n}}\bigr\|_{H^{-s}}&\lesssim \|f_1\|_{Y_{n,s}^1}\|f_2\|_{Y_{n,s}^1}\|g\|_{Y_{n,s}^1}^{\frac12}\|g\|_{H^1}^{\frac12}.\label{technical 3}
\end{align}
\end{lem}
\bpf
The estimates \eqref{technical 1} and \eqref{technical 2} are proved in \cite[Lemma 2.3]{2024arXiv240501246H}.
The estimate \eqref{technical 3} follows from a similar argument that we now give.

Recalling that \(Z_{n,s}(k) = \bigl\|\rho_k\tfrac{d\mu_n-1}{\sqrt{\epsilon_n}}\bigr\|_{H^{-s}}\), we first use the product estimate \eqref{product estimate} to bound
\begin{align*}
\bigl\|f_1f_2g\tfrac{d\mu_n-1}{\sqrt{\epsilon_n}}\bigr\|_{H^{-s}}&\lesssim\sum_{\substack{k_1,k_2,k_3,k_4\in \Z\\|k_j - k_\ell|\leq 1}}\|\rho_{k_1}f_1\times\rho_{k_2}f_2\times \rho_{k_3}g\|_{H^1}Z_{n,s}(k_4).
\end{align*}
Applying \eqref{GN}, Young's inequality, and \eqref{product estimate inf}, we have
\begin{align*}
&\|\rho_{k_1}f_1\times\rho_{k_2}f_2\times \rho_{k_3}g\|_{H^1}\\
&\qquad \lesssim \|\rho_{k_1}f_1\|_{L^\infty}\|\rho_{k_2}f_2\|_{L^\infty}\|\rho_{k_3}g\|_{H^1} + \|\rho_{k_1}f_1\|_{L^\infty}\|\rho_{k_2}f_2\|_{H^1}\|\rho_{k_3}g\|_{L^\infty}\\
&\qquad \quad + \|\rho_{k_1}f_1\|_{H^1}\|\rho_{k_2}f_2\|_{L^\infty}\|\rho_{k_3}g\|_{L^\infty}\\
&\qquad\lesssim \sum_{j=1}^2\|\rho_{k_j}f_j\|_{L^2}\|\rho_{k_{1-j}}f_{1-j}\|_{H^1}\|g\|_{H^1}\\
&\qquad\quad + \sum_{j=1}^2\|\rho_{k_j}f_j\|_{L^2}^{\frac12}\|\rho_{k_j}f_j\|_{H^1}^{\frac12}\|\rho_{k_3}g\|_{L^2}^{\frac12}\|\rho_{k_3}g\|_{H^1}^{\frac12}\|f_{1-j}\|_{H^1}.
\end{align*}

Recalling that if \(|k - k_4|\leq 1\) then \(Z_{n,s}(k_4)\lesssim \cN(k;Z_{n,s})\), we may combine these estimates with H\"older's inequality, \eqref{sq sum}, and \eqref{sq sum H1} to get
\begin{align*}
&\bigl\|f_1f_2g\tfrac{d\mu_n-1}{\sqrt{\epsilon_n}}\bigr\|_{H^{-s}}\\
&\qquad\lesssim \sum_{j=1}^2 \sum_{\substack{k_1,k_2\in \Z\\|k_1 - k_2|\leq 1}} \cN(k_1;Z_{n,s})\|\rho_{k_1}f_j\|_{L^2}\|\rho_{k_2}f_{1-j}\|_{H^1}\|g\|_{H^1}\\
&\qquad\quad + \sum_{j=1}^2\sum_{\substack{k_1,k_2\in \Z\\|k_1 - k_2|\leq 1}}\cN(k_1;Z_{n,s})^{\frac12}\|\rho_{k_1}f_j\|_{L^2}^{\frac12}\|\rho_{k_1}f_j\|_{H^1}^{\frac12}\\
&\qquad\quad\qquad\qquad\qquad\qquad\times\cN(k_2;Z_{n,s})^{\frac12}\|\rho_{k_2}g\|_{L^2}^{\frac12}\|\rho_{k_2}g\|_{H^1}^{\frac12}\|f_{1-j}\|_{H^1}\\
&\qquad\lesssim \sum_{j=1}^2 \left[\sum_{k\in \Z}|\cN(k;Z_{n,s})|^2\|\rho_k f_j\|_{L^2}^2\right]^{\frac12}\left[\sum_{k\in \Z}\|\rho_k f_{1-j}\|_{H^1}^2\right]^{\frac12}\|g\|_{H^1}\\
&\qquad\quad + \sum_{j=1}^2 \left[\sum_{k\in \Z}|\cN(k;Z_{n,s})|^2\|\rho_k f_j\|_{L^2}^2\right]^{\frac14}\left[\sum_{k\in \Z}\|\rho_k f_j\|_{H^1}^2\right]^{\frac14}\\
&\qquad\quad\qquad\qquad\times\left[\sum_{k\in \Z}|\cN(k;Z_{n,s})|^2\|\rho_k g\|_{L^2}^2\right]^{\frac14}\left[\sum_{k\in \Z}\|\rho_k g\|_{H^1}^2\right]^{\frac14}\|f_{1-j}\|_{H^1}\\
&\qquad \lesssim \sum_{j=1}^2\Bigl[ \|f_j\|_{L_{Z_{n,s}}^2}\|f_{1-j}\|_{H^1}\|g\|_{H^1} + \|f_j\|_{L_{Z_{n,s}}^2}^{\frac12}\|f_j\|_{H^1}^{\frac12}\|g\|_{L_{Z_{n,s}}^2}^{\frac12}\|g\|_{H^1}^{\frac12}\|f_{1-j}\|_{H^1}\Bigr].
\end{align*}
The estimate \eqref{technical 3} now follows from the definition \eqref{Yns1} of \(Y_{n,s}^1\).
\epf

We recall that the distribution of a random measure is uniquely determined by its Laplace functionals. The hypothesis \eqref{Levy exponential bound} on the L\'evy measure \(\Lambda\) appearing in \eqref{new Phi expression} suffices to ensure that \eqref{Laplace functional} holds for all bounded, compactly supported, Borel-measurable \(f\colon \R\to \C\) satisfying \(\Re f>-\frac a{\epsilon_n}\). This will be key to Corollaries~\ref{c:L2mu exp 1},~\ref{c:L2mu exp 2} and Proposition~\ref{p:Haar}, which lie at the heart of the proof of Theorems~\ref{t:homogenization} and~\ref{t:cvgce in distrib}. 

Let us start with the following generalization of \cite[Lemma 4.1]{2024arXiv240501246H}, which applies to any sequence of random variables \(Z\):
\begin{lem}\label{l:this one}
Let \(Z\colon \Z\to[0,\infty)\) be a sequence of non-negative random variables so that for all \(p\geq 1\) we have the estimate
\eq{uniform bound for Zs}{
\sup_{k\in \Z}\bbE |Z(k)|^{6p}\lesssim_p1.
}
Then for all \(p\geq 1\) and \(f\in L^2\) we have
\begin{align}
\sup_{k\in \Z}\bbE |\cN(k;Z)|^{2p}&\lesssim_p1,\label{exp Nk}\\
\bbE \|f\|_{L_Z^2}^{2p}&\lesssim_p\|f\|_{L^2}^{2p}.\label{LZ2 exp}
\end{align}
\end{lem}
\bpf
Proceeding as in \cite[Lemma 4.1]{2024arXiv240501246H}, we apply Young's inequality to bound
\[
|Z(\ell)|^2|k-\ell|^2\leq \tfrac 4{27}|Z(\ell)|^6+|k-\ell|^3,
\]
so that whenever \(k\neq \ell\) we have
\[
|Z(\ell)|^2 - |k-\ell|\lesssim\frac{|Z(\ell)|^6}{\<k-\ell\>^2}.
\]
Recalling the definition \eqref{defn: cnk}, for all \(p\geq 1\) we may then estimate
\[
|\cN(k;Z)|^{2p}\lesssim_p 1 + \sum_{\ell\in \Z}\frac{|Z(\ell)|^{6p}}{\<k-\ell\>^{2p}},
\]
and applying \eqref{uniform bound for Zs} yields
\[
\sup_{k\in \Z}\bbE |\cN(k;Z)|^{2p}\lesssim_p 1 + \sup_{\ell\in \Z}\bbE|Z(\ell)|^{6p}\lesssim_p1,
\]
which proves \eqref{exp Nk}.

Further, from \eqref{sq sum} and \eqref{sq sum H1} we have
\begin{align*}
\Bigl[\bbE \|f\|_{L_Z^2}^{2p}\Bigr]^{\frac1p} &\simeq \Biggl[\bbE\Biggl|\sum_{k\in \Z}|\cN(k;Z)|^2\|\rho_k f\|_{L^2}^2\Biggr|^p\Biggr]^{\frac1p}\notag\\
&\lesssim \sum_{k\in \Z} \Bigl[\bbE|\cN(k;Z)|^{2p}\Bigr]^{\frac1p}\|\rho_k f\|_{L^2}^2\lesssim_p\sum_{k\in \Z}\|\rho_kf\|_{L^2}^2\lesssim_p \|f\|_{L^2}^2,
\end{align*}
where we have used that the implicit constant in \eqref{sq sum} is independent of the sequence \(Z\). This completes the proof of \eqref{LZ2 exp}.
\epf

We now apply this lemma with \(Z = Z_n\) to obtain a uniform (in n) version of \cite[Lemma 4.1]{2024arXiv240501246H}:

\begin{cor}\label{c:L2mu exp 1}
Suppose that for all \(n\geq 1\) we have \(0<\epsilon_n\leq 1\). Taking \(Z_n\) to be defined as in \eqref{Xn1}, we almost surely have \(\cN(0;Z_n)<\infty\) for all \(n\geq 1\).

Further, if \(f\in H^1\) and \(p\geq 1\) then
\begin{align}
\sup_{n\geq 1}\bbE\|f\|_{X_n^1}^{2p}&\lesssim_p\|f\|_{H^1}^{2p}.\label{Xn1 exp}
\end{align}
In particular, we almost surely have \(f\in \bigcap_{n\geq 1}X_n^1\).
\end{cor}
\bpf 
We apply Lemma~\ref{l:this one} with \(Z = Z_n\). Employing \eqref{Laplace functional} with \(f = - s\bbo_{[k-\frac12,k+\frac12)}\) and recalling \eqref{new Phi expression}, we obtain the moment generating function
\[
\bbE \exp \bigl[sZ_n(k)\bigr] = \exp\bigl[- \tfrac1{\epsilon_n}\Phi(-\epsilon_n s) + s\bigr],
\]
provided \(|s|\ll1\) is sufficiently small. As \(0<\epsilon_n\leq 1\), if \(p\geq 1\) is an integer then we may use \eqref{Levy exponential bound} to obtain
\[
\bbE |Z_n(k)|^{6p} = \p_s^{6p} \bbE \exp \bigl[sZ_n(k)\bigr]\Bigr|_{s=0}\lesssim_p 1.
\]
By interpolation, we may extend this to all \(p\geq 1\) to obtain
\[
\sup_{n\geq 1}\sup_{k\in \Z}\bbE|Z_n(k)|^{6p} \lesssim_p1.
\]
The estimate \eqref{exp Nk} now ensures that almost surely \(\cN(0;Z_n)<\infty\) for all \(n\geq 1\). The estimate \eqref{Xn1 exp} follows from \eqref{LZ2 exp} and implies that \(f\in \bigcap_{n\geq 1}X_n^1\) almost surely.
\epf

For dyadic \(N\geq 1\) and integers \(k\in \Z\), we denote the dyadic interval
\[
I_{N,k} = \bigl[\tfrac kN,\tfrac{k+1}N\bigr).
\]
We define the associated Haar functions \(e_{N,k}\) by taking \(e_{1,k} = \bbo_{I_{1,k}}\) and
\[
e_{N,k} = \sqrt{\tfrac N2}\bigl(\bbo_{I_{N,2k}} - \bbo_{I_{N,2k+1}}\bigr)\qt{for \(N\geq 2\)}.
\]
We then have:
\begin{lem}\label{l:Haar series}
If \(f\in C_c(\R)\) then the partial sums
\eq{Haar series}{
\sum_{1\leq N\leq N_0}\sum_{k\in \Z}\<f,e_{N,k}\>e_{N,k}
}
converge uniformly to \(f\) as \(N_0\to\infty\).
\end{lem}
\bpf
Let
\[
E_N f = \sum_{k\in \Z} N\<f,\bbo_{I_{N,k}}\> \bbo_{I_{N,k}},
\]
where we note that as \(f\) has compact support the sum in \(k\) is finite.

For all \(N\geq 2\) we may decompose
\(
I_{\frac N2,k} = I_{N,2k}\sqcup I_{N,2k+1}
\),
and hence
\begin{align*}
E_Nf - E_{\frac N2}f &= \sum_{k\in \Z}\biggl\{\Bigl[N\<f,\bbo_{I_{N,2k}}\>-\tfrac N2\<f,\bbo_{I_{\frac N2,k}}\>\Bigr]\bbo_{N,2k}\\
&\qquad\qquad + \Bigl[N\<f,\bbo_{I_{N,2k+1}}\>-\tfrac N2\<f,\bbo_{I_{\frac N2,k}}\>\Bigr]\bbo_{N,2k+1}\biggr\}\\
&= \sum_{k\in \Z}\<f,e_{N,k}\>e_{N,k}.
\end{align*}
Consequently, for all \(N_0\geq 2\) we have
\[
E_{N_0} f = E_1 f + \sum_{2\leq N\leq N_0}\bigl[E_N f - E_{\frac N2}f\bigr] = \sum_{1\leq N\leq N_0}\sum_{k\in \Z} \<f,e_{N,k}\>e_{N,k}.
\]

As \(f\) is uniformly continuous on \(\R\), for all \(\eta>0\) we may find some dyadic \(N_1\geq 2\) so that whenever \(|x - y|\leq \frac1{N_1}\) we have
\[
|f(x) - f(y)|\leq \eta.
\]
For all dyadic \(N_0\geq N_1\) and all \(x\in \R\) we may then bound
\begin{align*}
\left|f(x) - \sum_{1\leq N\leq N_0}\sum_{k\in \Z}\<f,e_{N,k}\>e_{N,k}(x)\right| &= |f(x) - E_{N_0}f(x)| \\
&\leq \sum_{k\in \Z} N_0\left[\int_{I_{N_0,k}} |f(x) - f(y)|\,dy\right] \bbo_{I_{N_0,k}}(x)\\
&\leq\eta,
\end{align*}
which suffices to prove that \eqref{Haar series} converges to \(f\) uniformly on \(\R\).
\epf

For dyadic \(N_1,N_2\geq 1\) and any \(k_1,k_2\in \Z\), we compute that
\eq{Haar orthogonal}{
\<e_{N_1,k_1},e_{N_2,k_2}\> = \delta_{N_1,N_2}\delta_{k_1,k_2},
}
so \((e_{N,k})\) is an orthonormal subset of \(L^2\). In particular, Bessel's inequality ensures that for all \(\psi\in L^2\) the series
\(
\sum_{N\geq 1}\sum_{k\in \Z}\<\psi,e_{N,k}\>e_{N,k}
\)
converges absolutely in \(L^2\). Further, as \(C_c(\R)\) is dense in \(L^2\), Lemma~\ref{l:Haar series} shows that the set \((e_{N,k})\) is a Hilbert basis for \(L^2\).

We are now in a position to prove the following:
\begin{prop}\label{p:Haar}
Suppose that \(0<\epsilon_n\leq 1\). Then, for all \(s>\frac12\), \(\psi\in H^1\), and \(p\geq 1\) we have the estimate
\eq{Haar moments}{
\bbE \|\psi(d\mu_n-1)\|_{H^{-s}}^{2p}\lesssim_{p,s} \epsilon_n^p\|\psi\|_{H^1}^{2p}.
}
\end{prop}
\bpf
Let us fix the value of \(n\geq 1\) and denote the random variables
\[
X_{N,k} = \int_\R e_{N,k}\,d\mu_n - \int_\R e_{N,k}\,dx.
\]

As \(\bbE\mu_n\) is Lebesgue measure, we have
\eq{centered}{
\bbE X_{N,k} = 0.
}
If \(J\geq 2\) then, given dyadic \(N_1,\dots N_J\geq 1\), \(s_1,\dots,s_J\in \R\) so that \(|s_j|\ll1\) is sufficiently small, and integers \(k_1,\dots,k_J\in \Z\), we may use that \eqref{Laplace functional} extends to all bounded, compactly supported, Borel-measurable \(f\colon \R\to \C\) with \(\Re f>-a\) for \(a>0\) as in \eqref{Levy exponential bound}, to obtain
\[
\log \bbE \exp\Biggl[\sum_{j=1}^J s_j X_{N_j,k_j}\Biggr] = - \tfrac1{\epsilon_n}\int_\R \Phi\Biggl(- \epsilon_n\sum_{j=1}^J s_j e_{N_j,k_j}\Biggr)\,dx,
\]
where \(\Phi(z)\) is defined as in \eqref{new Phi expression}. Consequently, the joint cumulant
\begin{align}
\kappa\bigl(X_{N_1,k_1},\dots,X_{N_J,k_J}\bigr) &= \p_{s_1}\dots\p_{s_J} \log \bbE \exp\Biggl[\sum_{j=1}^J s_j X_{N_j,k_j}\Biggr]\Biggr|_{s_1=\dots = s_J=0} \notag\\
&=  (-\epsilon_n)^{J-1}\p_z^J\Phi\bigr|_{z=0}\int_\R \prod_{j=1}^J e_{N_j,k_j}\,dx.\label{cumulants}
\end{align}

By combining \eqref{centered} and \eqref{cumulants} with \eqref{Haar orthogonal} and our assumption that \(\Phi''(0) = -1\), we obtain the identity
\[
\bbE\bigl[X_{N_1,k_2}X_{N_2,k_2}\bigr] = \operatorname{cov}\bigl(X_{N_1,k_2},X_{N_2,k_2}\bigr) = \kappa\bigl(X_{N_1,k_1},X_{N_2,k_2}\bigr) = \epsilon_n\delta_{N_1,N_2}\delta_{k_1,k_2},
\]
so the random variables \(\frac1{\sqrt{\epsilon_n}}X_{N,k}\) are orthonormal.

To control the joint cumulants for larger values of \(J\), we assume (without loss of generality) that \(N_1\geq N_2\geq \dots\geq N_J\) and denote
\[
k_j' = k_j'(N_1,k_1,N_j) = \begin{cases}k_1&\text{if \(N_1 = 1\)},\\\bigl\lfloor\frac{2k_1}{N_1}\bigr\rfloor&\text{if \(1 = N_j<N_1\)},\\\bigl\lfloor\frac{N_j k_1}{N_1}\bigr\rfloor&\text{if \(1<N_j\leq N_1\)}.\end{cases}
\]
We then see that the supports of \(e_{N_1,k_1}\) and \(e_{N_j,k_j}\) are disjoint unless \(k_j = k_j'\). Further, if \(N_j<N_1\) then \(e_{N_j,k_j}\) is constant on the support of \(e_{N_1,k_1}\). In particular, if \(N_2<N_1\) then
\[
\int_\R e_{N_1,k
_1}\prod_{j=2}^J e_{N_j,k_j'}\,dx = 0.
\]
Combining these facts, we arrive at the inequality
\eq{higher cumulant bound}{
\bigl|\kappa\bigl(X_{N_1,k_1},\dots,X_{N_J,k_J}\bigr)\bigr| \lesssim_J \epsilon_n^{J-1} \delta_{N_1,N_2} \delta_{k_1,k_2}\prod_{j=3}^J \delta_{k_j,k_j'}N_j^{\frac12},
}
whenever \(J\geq 2\).

Next, we define the operator
\[
K = \<\p_x\>^{-s}\psi,
\]
and apply \cite[Theorem 4.1]{MR2154153} to estimate
\eq{K I2}{
\sum_{N\geq 1}\sum_{k\in \Z}\|K e_{N,k}\|_{L^2}^2 = \|K\|_{\I_2}^2\leq \|\<\xi\>^{-s}\|_{L^2}^2\|\psi\|_{L^2}^2\lesssim \|\psi\|_{L^2}^2.
}

For \(N\geq 2\) we compute that
\[
\int_{-\infty}^x e_{N,k}(y)\,dy = \sqrt{\tfrac N2}\bigl(x - \tfrac{2k}N\bigr) \bbo_{I_{N,2k}}(x) - \sqrt{\tfrac N2}\bigl(x - \tfrac{2k+2}N\bigr)\bbo_{I_{N,2k+1}}(x),
\]
and hence
\[
\|e_{N,k}\|_{H^{-1}}\lesssim \left\|\int_{-\infty}^x e_{N,k}(y)\,dy\right\|_{L^2} \lesssim \tfrac1N.
\]
Using also that \(\|e_{1,k}\|_{H^{-1}}\lesssim \|e_{1,k}\|_{L^2} = 1\), we may then interpolate to obtain
\[
\|e_{N,k}\|_{H^{-s}}\lesssim_s \|e_{N,k}\|_{H^{-1}}^s\|e_{N,k}\|_{L^2}^{1-s}\lesssim_s N^{-s},
\]
for all \(N\geq 1\). Applying \eqref{product estimate}, we may then bound
\eq{K Hs}{
\|Ke_{N,k}\|_{L^2} = \|\psi e_{N,k}\|_{H^{-s}}\lesssim_s \|\psi\|_{H^1}\|e_{N,k}\|_{H^{-s}}\lesssim_s N^{-s}\|\psi\|_{H^1}.
}

For a positive integer \(J\geq 1\), let us now denote
\[
S(J) = \sum_{\substack{(N_1,k_1),\dots,(N_J,k_J)\\N_1\geq \dots\geq N_J}} \bigl|\kappa\bigl(X_{N_1,k_1},\dots,X_{N_J,k_J}\bigr)\bigr| \prod_{j=1}^J \|Ke_{N_j,k_j}\|_{L^2}.
\]
From \eqref{centered}, we have \(S(1) = 0\). For \(J\geq 2\) we instead use \eqref{higher cumulant bound}, \eqref{K I2}, and \eqref{K Hs} to bound
\begin{align}
S(J) &\lesssim_J\epsilon_n^{J-1} \sum_{(N_1,k_1)} \|Ke_{N_1,k_1}\|_{L^2}^2\sum_{N_3,\dots,N_J} \prod_{j=3}^J N_j^{\frac12}\|Ke_{N_j,k_j'}\|_{L^2}\notag\\
&\lesssim_{J,s} \epsilon_n^{J-1} \|\psi\|_{H^1}^J\sum_{N_3,\dots,N_J} \prod_{j=3}^J N_j^{\frac12-s}\notag\\
&\lesssim_{J,s}  \epsilon_n^{J-1}\|\psi\|_{H^1}^J.\label{Sj bound}
\end{align}

Given an integer \(p\geq 1\), the moment-cumulant formula (see, e.g., \cite{MR123345}) tells us that
\[
\bbE\prod_{j=1}^{2p}X_{N_j,k_j} = \sum_{P\in \mc P_{2p}}\prod_{B\in P}\kappa\bigl((X_{N_j,k_j})_{j\in B}\bigr),
\]
where \(\mc P_{2p}\) is the set of all partitions of the set \(\{1,\dots,2p\}\). 
Combining this with \eqref{Sj bound} and using that \(\epsilon_n\leq 1\), we obtain the estimate
\begin{align}
\sum_{(N_1,k_1),\dots,(N_{2p},k_{2p})} \Biggl|\bbE\prod_{j=1}^{2p}X_{N_j,k_j}\Biggr| \prod_{\ell=1}^{2p} \|Ke_{N_\ell,k_\ell}\|_{L^2}&\lesssim_p \sum_{J_1 + J_2 + \dots + J_L = 2p}\prod_{\ell=1}^L S(J_\ell)\label{big bound}\\
&\lesssim_{p,s} \epsilon_n^p \|\psi\|_{H^1}^{2p},\notag
\end{align}
where the sum in the middle expression is interpreted to be over all positive integers \(L\geq 1\) and \(J_\ell\geq 1\) so that \(J_1+\dots+J_L = 2p\).

Given a finite set
\(
F\subseteq \bigl\{(N,k):N\geq1\text{ is dyadic and }k\in \Z\bigr\}
\),
we compute that
\begin{align*}
&\bbE\Biggl\|\psi\sum_{(N,k)\in F}X_{N,k}e_{N,k}\Biggr\|_{H^{-s}}^{2p}\\
&\qquad = \sum_{(N_1,k_1),\dots,(N_{2p},k_{2p})\in F} \bbE \prod_{j=1}^{2p}X_{N_j,k_j} \prod_{\ell=1}^p \<Ke_{N_{2\ell-1},k_{2\ell-1}},Ke_{N_{2\ell},k_{2\ell}}\>\\
&\qquad\leq \sum_{(N_1,k_1),\dots,(N_{2p},k_{2p})\in F} \Biggl|\bbE\prod_{j=1}^{2p}X_{N_j,k_j}\Biggr| \prod_{\ell=1}^{2p} \|Ke_{N_\ell,k_\ell}\|_{L^2}.
\end{align*}
Combining this with the (absolute) convergence of the series \eqref{big bound}, we see that the partial sums of any enumeration of the series \(\psi \sum_{N\geq 1}\sum_{k\in \Z}X_{N,k}e_{N,k}\) are Cauchy in \(L^{2p}(d\bbP;H^{-s})\). Consequently, this series converges unconditionally in \(L^{2p}(d\bbP;H^{-s})\) and the corresponding limit satisfies
\eq{Haar moments basically}{
\bbE\left\|\psi\sum_{N\geq 1}\sum_{k\in \Z} X_{N,k}e_{N,k}\right\|_{H^{-s}}^{2p} \lesssim_{s,p} \epsilon_n^p\|\psi\|_{H^1}^{2p}.
}
As this holds for all integers \(p\geq 1\), we may extend to all real numbers \(p\geq 1\) by interpolation.

Let us now assume that \(\psi\in C_c^\infty(\R)\).
If \(\phi\in H^s\) then \(\phi\psi\in C_c(\R)\) and Lemma~\ref{l:Haar series} gives us
\[
\int_\R \phi\psi\,d\mu_n - \int_\R \phi\psi\,dx = \lim_{N_0\to\infty} \int_\R \phi\psi \sum_{1\leq N\leq N_0}\sum_{k\in \Z}X_{N,k}e_{N,k}\,dx,
\]
almost surely. Combining this with the unconditional convergence of the series \(\psi \sum_{N\geq 1}\sum_{k\in \Z}X_{N,k}e_{N,k}\) in \(L^{2p}(d\bbP;H^{-s})\), we see that almost surely
\[
\psi\sum_{N\geq 1}\sum_{k\in \Z}X_{N,k}e_{N,k} = \psi(d\mu_n-1),
\]
in \(H^{-s}\) and \eqref{Haar moments} follows from \eqref{Haar moments basically}. Finally, the extension of \eqref{Haar moments} to all \(\psi\in H^1\) follows from the density of \(C_c^\infty\) in \(H^1\).
\epf

Combining the estimate \eqref{Haar moments} with the proof of Corollary~\ref{c:L2mu exp 1} gives us:

\begin{cor}\label{c:L2mu exp 2}
Suppose that for all \(n\geq 1\) we have \(0<\epsilon_n\leq 1\). Taking \(\frac12<s\leq 1\) and \(Z_{n,s}\) to be defined as in \eqref{Yns1}, we almost surely have \(\cN(0;Z_{n,s})<\infty\) for all \(n\geq 1\).

Further, if \(f\in H^1\) and \(p\geq 1\) then
\begin{align}
\sup_{n\geq 1}\bbE\|f\|_{Y_{n,s}^1}^{2p}&\lesssim_{p,s}\|f\|_{H^1}^{2p},\label{Yns1 exp}
\end{align}
and hence we almost surely have \(f\in \bigcap_{n\geq 1}Y_{n,s}^1\).
\end{cor}
\bpf
Applying \eqref{Haar moments}, we have
\eq{pug}{
\sup_{n\geq 1}\sup_{k\in \Z}\bbE |Z_{n,s}(k)|^{6p}\lesssim_{p,s} \sup_{k\in \Z}\|\rho_k\|_{H^1}^{6p}\lesssim_{p,s}1.
}
Arguing as in Corollary~\ref{c:L2mu exp 1}, we combine \eqref{pug} with \eqref{exp Nk} to show that \(\cN(0;Z_{n,s})<\infty\) for all \(n\geq 1\) almost surely. The estimate \eqref{Yns1 exp} follows from \eqref{Xn1 exp} and combining \eqref{pug} with \eqref{LZ2 exp}. The fact that \(f\in \bigcap_{n\geq 1}Y_{n,s}^1\) almost surely is a consequence of \eqref{Yns1 exp}.
\epf

We conclude this section by applying a similar argument to Proposition~\ref{p:Haar} (see also \cite[Equation (2.12)]{MR4145790}) to obtain the following:
\begin{lem}\label{l:expect white noise}
If \(\xi\) is white-noise-distributed then for all \(s>\frac12\) and \(\psi\in L^2\) we have the estimate
\eq{expected white noise}{
\bbE \|\psi\xi\|_{H^{-s}}^2\lesssim \|\psi\|_{L^2}^2.
}
\end{lem}
\bpf
We recall (see, e.g., \cite[Theorem 2.7]{MR4145790}) that if \((e_j)\) is a real-valued orthonormal basis of \(L^2\) and \((X_j)\) is an i.i.d. sequence of standard normal random variables then the series
\[
\xi = \sum_{j=1}^\infty X_je_j
\]
converges almost surely in the sense of distributions to white noise.

Taking \(K = \<\p_x\>^{-s}\psi\) as in Proposition~\ref{p:Haar}, for any finite set \(F\subseteq \Z^+\) we have
\begin{align*}
\bbE\biggl\|\psi \sum_{j\in F} X_j e_j\biggr\|_{H^{-s}}^2 = \sum_{j_1,j_2\in F}\bbE X_{j_1}X_{j_2}\<Ke_{j_1}, K e_{j_2}\> = \sum_{j\in F}\|Ke_j\|_{L^2}^2.
\end{align*}
Using that
\[
\sum_{j=1}^\infty\|Ke_j\|_{L^2}^2 = \|K\|_{\I_2}^2\lesssim \|\psi\|_{L^2}^2,
\]
an identical argument to Proposition~\ref{p:Haar} now shows that if \(\psi\in \Test\) then \(\psi\sum_{j=1}^\infty X_je_j\) converges unconditionally in \(L^2(d\bbP;H^{-s})\) to \(\psi\xi\) and satisfies
\eq{expected white noise almost}{
\bbE\biggl\|\psi\sum_{j=1}^\infty X_je_j\biggr\|_{H^{-s}}^2 \lesssim \|\psi\|_{L^2}^2.
}
The estimate \eqref{expected white noise} then follows from \eqref{expected white noise almost} and the density of \(\Test\) in \(L^2\).
\epf

\section{Homogenization}\label{s:homogenization}

We start this section by recording some properties of the solution \(\psi_n\) of \eqref{NLSn} obtained in \cite{2024arXiv240501246H}:
\begin{lem}\label{l:ap psin}
Suppose that \(\cN(0;Z_n)<\infty\) and \(\psi_0\in X_n^1\). Then there exists a unique solution \(\psi_n\in C(\R;X_n^1)\) of the Duhamel formulation,
\eq{Duhamel NLSn}{
\psi_n(t) = e^{it\p_x^2}\psi_0 - 2i\int_0^t e^{i(t-\tau)\p_x^2}\bigl[|\psi_n(\tau)|^2\psi_n(\tau)\,d\mu_n\bigr]\,d\tau,
}
of \eqref{NLSn} with initial data \(\psi_n(0) = \psi_0\). This solution conserves both the mass and energy,
\eq{mass and energy n}{
M\bigl[\psi_n\bigr] = \int_\R|\psi_n|^2\,dx \qtq{and}E_n\bigl[\psi_n\bigr] = \tfrac12\int_\R |\p_x\psi_n|^2\,dx + \tfrac12\int_\R |\psi_n|^4\,d\mu_n,
}
and for all \(T>0\) satisfies the estimates
\begin{align}
\sup_{t\in \R}\|\psi_n(t)\|_{H^1}^2 &\lesssim \bigl[1 + \|\psi_0\|_{X_n^1}^2\bigr]\|\psi_0\|_{H^1}^2,\label{psin H1}\\
\|\psi_n\|_{C_TX_n^1}^2 &\lesssim \<T\>\bigl[1 + \|\psi_0\|_{X_n^1}^2\bigr]\|\psi_0\|_{X_n^1}^2.\label{ap for psin}
\end{align}
\end{lem}
\bpf
The existence and uniqueness of a solution follows from \cite[Proposition 3.5]{2024arXiv240501246H}. The estimates \eqref{psin H1} and \eqref{ap for psin} follow from combining \cite[Lemma 3.3]{2024arXiv240501246H} with \cite[Proposition 3.5]{2024arXiv240501246H}.
\epf

We will also require a corresponding result for the solution \(\psi\) of \eqref{NLS}. This follows from the standard well-posedness theory for \eqref{NLS} and arguments appearing in \cite[Lemma 3.3]{2024arXiv240501246H}:
\begin{lem}\label{l:ap psi}
If \(\psi_0\in H^1\) then there exists a unique solution \(\psi\in C(\R;H^1)\) of the Duhamel formulation,
\eq{Duhamel NLS}{
\psi(t) = e^{it\p_x^2}\psi_0 - 2i\int_0^t e^{i(t-\tau)\p_x^2}\bigl[|\psi(\tau)|^2\psi(\tau)\bigr]\,d\tau,
}
of \eqref{NLS} with initial data \(\psi(0) = \psi_0\). This solution conserves both the mass and energy,
\eq{mass and energy}{
M[\psi] = \int_\R |\psi|^2\,dx \qtq{and}E[\psi] = \tfrac12\int_\R |\p_x\psi|^2 + |\psi|^4\,dx,
}
satisfies the estimate
\begin{align}
\sup_{t\in \R}\|\psi(t)\|_{H^1}^2 &\lesssim \bigl[1 + \|\psi_0\|_{L^2}^2\bigr]\|\psi_0\|_{H^1}^2,\label{psi H1}
\end{align}
and for all \(T>0\) we have
\begin{align}
\lim_{R\to\infty}\|(1 - \chi_R)\psi\|_{C_TH^1} &=0 = \lim_{N\to\infty}\|P_{>N}\psi\|_{C_TH^1}.\label{tight in H1}
\end{align}

Further, if \(\frac12<s\leq 1\), both \(\cN(0;Z_n)<\infty\) and \(\cN(0;Z_{n,s})<\infty\), and \(\psi_0\in Y_{n,s}^1\) then \(\psi\in C(\R;Y_{n,s}^1)\) and for all \(T>0\) we have the estimate
\begin{align}
\|\psi\|_{C_TY_{n,s}^1}^2 &\lesssim \<T\>\bigl[1 + \|\psi_0\|_{L^2}^2\bigr]\|\psi_0\|_{Y_{n,s}^1}^2.\label{ap for psi}
\end{align}
\end{lem}
\bpf
Global well-posedness as well as conservation of mass and energy follow from classical arguments; see, e.g., \cite{MR2002047}. We may then combine the conservation of mass and energy with \eqref{GN} to obtain \eqref{psi H1}: For all \(t\in \R\) we have
\[
\|\psi(t)\|_{H^1}^2 \leq M\bigl[\psi(t)\bigr] + 2E\bigl[\psi(t)\bigr] = M\bigl[\psi_0\bigr] + 2E\bigl[\psi_0\bigr]\lesssim \RHS{psi H1}.
\]
The identity \eqref{tight in H1} follows from applying Lemma~\ref{l:compactness helper} with \(\mc F = \{\psi\}\).

To prove that \(\psi\in C(\R;Y_{n,s}^1)\) and obtain the estimate \eqref{ap for psi}, we essentially follow the argument of \cite[Lemma 3.3]{2024arXiv240501246H}.

We first note that if \(\psi\in C(\R;H^1)\) satisfies \eqref{Duhamel NLS} then \(\psi\in C^1(\R;H^{-1})\) satisfies the equation \eqref{NLS} in \(H^{-1}\) for all \(t\in \R\).
Let \(Z\colon \Z\to[0,\infty)\) be a sequence of non-negative numbers and suppose that both \(\cN(0;Z)<\infty\) and \(\psi_0\in L_Z^2\). Taking \(R\geq 1\) and \(\chi_R\) to be defined as in Section~\ref{s:prelim}, we then compute that
\begin{align*}
&\tfrac d{dt}\|\chi_R \psi(t)\|_{L_Z^2}^2\\
&\qquad= -2\Im\int_\R \p_x^2\psi(t,x)\overline{\psi(t,x)} |\chi_R(x)|^2 \omega(x;Z)\,dx\\
&\qquad = 2\Im \int_\R \p_x\psi(t,x)\overline{\psi(t,x)}\bigl[2\chi_R(x)\chi_R'(x)\omega(x;Z) + |\chi_R(x)|^2 \p_x\omega(x;Z)\bigr]\,dx.
\end{align*}

Recalling that the map \(x\mapsto \omega(x;Z)\) is \(1\)-Lipschitz, as well as the estimate \eqref{Lipschitz bound}, we may use the conservation of mass and \eqref{psi H1} to bound
\begin{align*}
\Bigl|\tfrac d{dt}\|\chi_R \psi(t)\|_{L_Z^2}^2\Bigr|&\lesssim \Bigl[1 + \tfrac{|\cN(0;Z)|^2} R\Bigr]\|\psi(t)\|_{H^1}\|\chi_R\psi(t)\|_{L^2}\\
&\lesssim \Bigl[1 + \tfrac{|\cN(0;Z)|^2} R\Bigr]\bigl[1 + \|\psi_0\|_{L^2}^2\bigr]^{\frac12}\|\psi_0\|_{H^1}\|\psi_0\|_{L^2}.
\end{align*}
Integrating this expression, for all \(t\in \R\) we have
\[
\|\chi_R\psi(t)\|_{L_Z^2}^2 \lesssim \|\psi_0\|_{L_Z^2}^2 + |t| \Bigl[1 + \tfrac{|\cN(0;Z)|^2} R\Bigr]\bigl[1 + \|\psi_0\|_{L^2}^2\bigr]^{\frac12}\|\psi_0\|_{H^1}\|\psi_0\|_{L^2}.
\]
Taking \(R\to\infty\) using Fatou's lemma, we arrive at the estimate
\eq{LTinfLZ2}{
\sup_{|t|\leq T}\|\psi(t)\|_{L_Z^2}^2\lesssim \|\psi_0\|_{L_Z^2}^2 + T\bigl[1 + \|\psi_0\|_{L^2}^2\bigr]^{\frac12}\|\psi_0\|_{H^1}\|\psi_0\|_{L^2},
}
which we note is independent of the value of \(\cN(0;Z)\).

Repeating this argument, for \(|t|\leq T\) we have the estimate
\[
\Bigl|\tfrac d{dt}\|(1 - \chi_R)\psi(t)\|_{L_Z^2}^2 \Bigr|\lesssim \Bigl[1 + \tfrac{|\cN(0;Z)|^2} R\Bigr]\bigl[1 + \|\psi_0\|_{L^2}^2\bigr]^{\frac12}\|\psi_0\|_{H^1}\|(1 - \chi_R)\psi\|_{C_TL^2},
\]
and hence
\begin{align*}
&\sup_{|t|\leq T}\|(1 - \chi_R)\psi(t)\|_{L_Z^2}^2\\
&\qquad\lesssim \|(1 - \chi_R)\psi_0\|_{L_Z^2}^2 + T\Bigl[1 + \tfrac{|\cN(0;Z)|^2} R\Bigr]\bigl[1 + \|\psi_0\|_{L^2}^2\bigr]^{\frac12}\|\psi_0\|_{H^1}\|(1 - \chi_R)\psi\|_{C_TL^2},
\end{align*}
so, applying the Dominated Convergence Theorem to the first summand and \eqref{tight in H1} to the second, yields
\[
\lim_{R\to\infty}\sup_{|t|\leq T}\|(1 - \chi_R)\psi(t)\|_{L_Z^2} = 0.
\]
We may then combine this with \eqref{Lipschitz bound} to show that for all \(t,\tau\in [-T,T]\) we have
\begin{align*}
\|\psi(t) - \psi(\tau)\|_{L_Z^2} &\leq \bigl\|\chi_R\bigl[\psi(t) - \psi(\tau)\bigr]\bigr\|_{L_Z^2} + 2\sup_{|t|\leq T}\|(1 - \chi_R)\psi(t)\|_{L_Z^2}\\
&\lesssim \bigl[|\cN(0;Z)|^2 + R\bigr]\|\psi(t) - \psi(\tau)\|_{L^2} + 2\sup_{|t|\leq T}\|(1 - \chi_R)\psi(t)\|_{L_Z^2}.
\end{align*}
By first taking \(\tau\to t\) and then \(R\to \infty\), we see that \(\psi\in C(\R;L_Z^2)\).

Finally, we apply these results first with \(Z = Z_n\) and then with \(Z = Z_{n,s}\) to prove that whenever \(\cN(0;Z_n)<\infty\), \(\cN(0;Z_{n,s})<\infty\), and \(\psi_0\in Y_{n,s}^1\) we have \(\psi\in C(\R;Y_{n,s}^1)\), with the estimate \eqref{ap for psi} following from \eqref{psi H1} and \eqref{LTinfLZ2}.
\epf

We now apply our a priori estimates to obtain the following proposition, the proof of which relies on a similar circle of estimates to \cite[Proposition 3.5]{2024arXiv240501246H}:
\begin{prop}\label{p:the beast}
Suppose that \(0<\epsilon_n\leq 1\), \(\cN(0;Z_n)<\infty\), \(\cN(0;Z_{n,1})<\infty\), \(\psi_0\in Y_{n,1}^1\), and for some \(T,A>0\) we have 
\eq{A defn}{
\|\psi_n\|_{C_TX_n^1}^2 + \|\psi\|_{C_TY_{n,1}^1}^2\leq A^2.
}
Then there exists a constant \(C>0\) so that
\begin{align}
\|\psi_n - \psi\|_{C_TH^{-1}}&\lesssim e^{CA^8T^2}A^3T\bigl[1 + A^2T\bigr]\epsilon_n^{\frac18},\label{H-1}\\
\sup_{|t|\leq T}\Bigl|\|\psi_n(t)\|_{H^1}^2 - \|\psi(t)\|_{H^1}^2\Bigr|&\lesssim e^{CA^8T^2}A^4\bigl[1 + A^2T\bigr]\epsilon_n^{\frac18}.\label{energy diff}
\end{align}
\end{prop}
\bpf
We start by writing the difference
\begin{align*}
\psi_n(t) - \psi(t) &= -2i\bigl[\cI_1(t)+\cI_2(t)\bigr],
\end{align*}
where
\begin{align*}
\cI_1(t)&=\int_0^t e^{i(t-\tau)\p_x^2}\Bigl[\bigl(|\psi_n(\tau)|^2\psi_n(\tau) - |\psi(\tau)|^2\psi(\tau)\bigr)\,d\mu_n\Bigr]\,d\tau,\\
\cI_2(t)&=\int_0^t e^{i(t-\tau)\p_x^2}\Bigl[|\psi(\tau)|^2\psi(\tau)\,(d\mu_n-1)\Bigr]\,d\tau.
\end{align*}

Taking a test function \(f\in L^1\), we recall that for \(0<\tau<t\) we have the dispersive estimate
\eq{dispersive}{
\bigl\|e^{-i(t-\tau)\p_x^2}f\bigr\|_{L^\infty}\lesssim \tfrac1{\sqrt{t - \tau}}\|f\|_{L^1}.
}
Taking \(0<t\leq T\), we then apply \eqref{technical 1} followed by \eqref{dispersive} and H\"older's inequality to bound
\begin{align*}
|\<\cI_1(t),f\>|&\lesssim \int_0^t \int_\R\bigl|e^{-i(t-\tau)\p_x^2}f\bigr| \bigl[|\psi_n(\tau)|^2 + |\psi(\tau)|^2\bigr]|\psi_n(\tau) - \psi(\tau)|\,d\mu_n\,d\tau\\
&\lesssim A^2\|f\|_{L^1}\int_0^t \tfrac1{\sqrt{t-\tau}} \|\psi_n(\tau) - \psi(\tau)\|_{L^\infty}\,d\tau\\
&\lesssim A^2\|f\|_{L^1}t^{\frac14}\left[\int_0^t \|\psi_n(\tau) - \psi(\tau)\|_{L^\infty}^4\,d\tau\right]^{\frac14}.
\end{align*}
By duality, this yields the estimate
\eq{I1 inf}{
\|\cI_1(t)\|_{L^\infty}^4 \lesssim  A^8t\int_0^t \|\psi_n(\tau) - \psi(\tau)\|_{L^\infty}^4\,d\tau.
}

Applying \eqref{technical 2}, we may instead estimate
\eq{I1 H-1}{
\|\cI_1\|_{C_TH^{-1}}\lesssim A^2T\|\psi_n - \psi\|_{C_TL^\infty},
}
and from \eqref{technical 3} we also have
\eq{I2 H-1}{
\|\cI_2\|_{C_TH^{-1}}\lesssim A^3T\sqrt{\epsilon_n}.
}

We now combine these estimates to obtain \eqref{H-1}. For dyadic \(N\geq 1\) and times \(0<t\leq T\), we use Bernstein's inequality, followed by \eqref{I1 inf} and \eqref{I2 H-1} to obtain
\begin{align*}
\|\psi_n(t) - \psi(t)\|_{L^\infty}^4 \!&\lesssim \|P_{\leq N}\cI_1(t)\|_{L^\infty}^4 \!+ \|P_{\leq N}\cI_2(t)\|_{L^\infty}^4 + \bigl\|P_{>N}\bigl[\psi_n(t) - \psi(t)\bigr]\bigr\|_{C_TL^\infty}^4\\
&\lesssim \|\cI_1(t)\|_{L^\infty}^4 + N^6\|\cI_2\|_{C_TH^{-1}}^4 +  N^{-2}\bigl[\|\psi_n\|_{C_TH^1}^4 + \|\psi\|_{C_TH^1}^4\bigr]\\
&\lesssim A^8t\int_0^t \|\psi_n(\tau) - \psi(\tau)\|_{L^\infty}^4\,d\tau + N^6A^{12}T^4\epsilon_n^2 + N^{-2}A^4.
\end{align*}
Noting that the left hand side of this inequality is simply \(\frac d{dt}\int_0^t \||\psi_n(\tau) - \psi(\tau)\|_{L^\infty}^4\,d\tau\), we may apply Gronwall's inequality to bound
\begin{align*}
\int_0^T \|\psi_n(t) - \psi(t)\|_{L^\infty}^4\,dt\lesssim  e^{CA^8T^2}T \Bigl[N^6 A^{12}T^4\epsilon_n^2 + N^{-2}A^4\Bigr],
\end{align*}
for some constant \(C>0\). After applying a similar argument for negative times and possibly increasing the size of the constant \(C\), we arrive at the estimate
\begin{align*}
\|\psi_n - \psi\|_{C_T L^\infty} &\lesssim A^2T^{\frac14}\|\psi_n - \psi\|_{L_T^4L^\infty} + N^{\frac32}A^3T\sqrt{\epsilon_n} + N^{-\frac12}A\\
&\lesssim e^{CA^8T^2}\Bigl[N^{\frac32}A^3T\sqrt{\epsilon_n} + N^{-\frac12}A\Bigr].
\end{align*}
As this estimate is uniform in \(N\), we may take \(N\simeq \epsilon_n^{-\frac14}\) to arrive at
\eq{Gronwall inf}{
\|\psi_n - \psi\|_{C_T L^\infty}\lesssim e^{CA^8T^2}A\bigl[1 + A^2T\bigr]\epsilon_n^{\frac18}.
}
Finally, by combining \eqref{I1 H-1}, \eqref{I2 H-1}, and \eqref{Gronwall inf} we get
\begin{align*}
\|\psi_n - \psi\|_{C_TH^{-1}} &\lesssim \|\cI_1\|_{C_TH^{-1}} + \|\cI_2\|_{C_TH^{-1}}\\
&\lesssim A^2T\|\psi_n - \psi\|_{C_TL^\infty} + A^3T\sqrt{\epsilon_n}\lesssim\RHS{H-1}.
\end{align*}

It remains to prove \eqref{energy diff}. Taking \(t\in [-T,T]\) and \(R\geq 1\), we use conservation of mass and energy to bound
\begin{align*}
&\Bigl|\|\psi_n(t)\|_{H^1}^2 - \|\psi(t)\|_{H^1}^2\Bigr|\\
&\lesssim \Bigl|M[\psi_n(t)] - M[\psi(t)]\Bigr| + \Bigl|E_n[\psi_n(t)] - E[\psi(t)]\Bigr| + \left|\int_\R |\psi_n(t)|^4\,d\mu_n - \int_\R |\psi(t)|^4\,dx\right|\\
&\lesssim \left|\int_\R |\psi_0|^4\,d\mu_n - \int_\R |\psi_0|^4\,dx\right| + \left|\int_\R |\psi_n(t)|^4\,d\mu_n - \int_\R |\psi(t)|^4\,dx\right|\\
&\lesssim \sup_{|t|\leq T}\left|\int_\R |\psi_n(t)|^4 - |\psi(t)|^4\,d\mu_n\right| + \sup_{|t|\leq T}\left|\int_\R |\psi(t)|^4 \,d\mu_n - \int_\R|\psi(t)|^4\,dx\right|.
\end{align*}

The first term may be estimated using \eqref{technical 1} and \eqref{GN} followed by \eqref{Gronwall inf} to obtain
\begin{align*}
\sup_{|t|\leq T}\left|\int_\R |\psi_n(t)|^4-|\psi(t)|^4\,d\mu_n\right|&\lesssim A^3\|\psi_n-\psi\|_{C_TL^\infty}\lesssim e^{CA^8T^2}A^4\bigl[1 + A^2T\bigr]\epsilon_n^{\frac18}.
\end{align*}

For the second term, we first note that from the definition \eqref{Yns1} and the estimates \eqref{GN} and \eqref{product estimate} we may bound
\[
\|fg\|_{Y_{n,1}^1}\lesssim \|f\|_{Y_{n,1}^1}\|g\|_{Y_{n,1}^1}.
\]
Consequently, we may apply \eqref{technical 3} to obtain
\begin{align*}
\sup_{|t|\leq T}\left|\int_\R |\psi(t)|^4 \,d\mu_n - \int_\R|\psi(t)|^4\,dx\right| &\lesssim A^4\sqrt{\epsilon_n}.
\end{align*}
Combining these estimates yields \eqref{energy diff}.
\epf

We conclude this section with the
\bpf[Proof of Theorem~\ref{t:homogenization}] Let \(T>0\) be fixed and, after possibly passing to a tail of the sequence, suppose that \(0<\epsilon_n\leq 1\) for all \(n\geq 1\).

We first note that by applying \eqref{ap for psin} and \eqref{ap for psi}, followed by \eqref{Xn1 exp} and \eqref{Yns1 exp}, for all \(p\geq 1\) we have
\begin{align}
\sup_{n\geq 1}\bbE \Bigl[\|\psi_n\|_{C_TX_n^1}^{2p} + \|\psi\|_{C_TY_{n,1}^1}^{2p}\Bigr] &\lesssim_{T,p}\bbE\Bigl[\bigl[1 + \|\psi_0\|_{X_n^1}^{2p}\bigr]\|\psi_0\|_{Y_{n,1}^1}^{2p}\Bigr]\notag\\
&\lesssim_{T,p}\bigl[1 + \|\psi_0\|_{H^1}^{2p}\bigr]\|\psi_0\|_{H^1}^{2p}.\label{uniform bounds for expectation in Xn1}
\end{align}

Given a constant \(A>0\), we now take \(\Omega_{n,A}\) to be the event that \(\cN(0;Z_n)<\infty\), \(\cN(0;Z_{n,1})<\infty\), \(\psi_0\in Y_{n,1}^1\), and \eqref{A defn} is true. By Markov's inequality and \eqref{uniform bounds for expectation in Xn1} we then have
\begin{align}
\bbP\bigl[\Omega_{n,A}^c\bigr]&\leq A^{-2}\bbE \Bigl[\|\psi_n\|_{C_TX_n^1}^2 + \|\psi\|_{C_TY_{n,1}^1}^2\Bigr]\lesssim_T A^{-2}\bigl[1 + \|\psi_0\|_{H^1}^2\bigr]\|\psi_0\|_{H^1}^2.\label{omega n A c}
\end{align}

Applying \eqref{H-1}, we get
\begin{align*}
\bbE\Bigl[\|\psi_n - \psi\|_{C_TH^{-1}}^{2p}\mathbbm 1_{\Omega_{n,A}}\Bigr] &\lesssim_{A,T,p}\epsilon_n^{\frac p4}.
\end{align*}
In the complement of this event, conservation of mass and \eqref{omega n A c} give us
\begin{align*}
\bbE\Bigl[ \|\psi_n - \psi\|_{C_TH^{-1}}^{2p}\mathbbm 1_{\Omega_{n,A}^c}\Bigr] &\lesssim_p\|\psi_0\|_{L^2}^{2p}\bbP\bigl[\Omega_{n,A}^c\bigr]\\
&\lesssim_{T,p} A^{-2}\bigl[1 + \|\psi_0\|_{H^1}^2\bigr]\|\psi_0\|_{H^1}^2\|\psi_0\|_{L^2}^{2p}.
\end{align*}
Combining these estimates, we arrive at
\begin{align*}
\bbE\|\psi_n - \psi\|_{C_TH^{-1}}^{2p} &\leq \bbE\Bigl[\|\psi_n - \psi\|_{C_TH^{-1}}^{2p}\mathbbm 1_{\Omega_{n,A}}\Bigr] + \bbE\Bigl[ \|\psi_n - \psi\|_{C_TH^{-1}}^{2p}\mathbbm 1_{\Omega_{n,A}^c}\Bigr]\\
&\lesssim_{T,p} C(A,T,p) \epsilon_n^{\frac p4} + A^{-2}\bigl[1 + \|\psi_0\|_{H^1}^2\bigr]\|\psi_0\|_{H^1}^2\|\psi_0\|_{L^2}^{2p},
\end{align*}
for some constant \(C(A,T,p)>0\). Taking \(n\to\infty\) and then \(A\to\infty\) yields the identity
\eq{H-1 convergence}{
\lim_{n\to\infty}\bbE\|\psi_n - \psi\|_{C_TH^{-1}}^{2p} = 0.
}

To upgrade this to convergence in \(H^1\), we will use the Radon--Riesz Theorem. Denoting the \(H^1\)-inner product by \(\<f,g\>_{H^1}\), we first bound
\begin{align*}
&\bbE\sup_{|t|\leq T}\bigl|\<\psi_n - \psi,\psi\>_{H^1}\bigr|^{2p}\\
&\quad \lesssim_p \bbE\sup_{|t|\leq T}\bigl|\<P_{\leq N}(\psi_n - \psi),\psi\>_{H^1}\bigr|^{2p} + \bbE\sup_{|t|\leq T}\bigl|\<\psi_n - \psi,P_{>N}\psi\>_{H^1}\bigr|^{2p}\\
&\quad\lesssim_p N^{2p}\bbE\|\psi_n - \psi\|_{C_TH^{-1}}^{2p}\|\psi\|_{C_TH^1}^{2p} + \Bigl[\bbE\|\psi_n\|_{C_TH^1}^{2p} + \|\psi\|_{C_TH^1}^{2p}\Bigr]\|P_{>N}\psi\|_{C_TH^1}^{2p}.
\end{align*}
Employing \eqref{psi H1} and \eqref{H-1 convergence} for the first term, and \eqref{psi H1}, \eqref{uniform bounds for expectation in Xn1}, and \eqref{tight in H1} for the second, we may take \(n\to\infty\) and then \(N\to \infty\) to get
\eq{duality convergence}{
\lim_{n\to\infty}\bbE\sup_{|t|\leq T}\bigl|\<\psi_n - \psi,\psi\>_{H^1}\bigr|^{2p} = 0.
}

Turning to the difference of the norms, from \eqref{energy diff} we have
\begin{align*}
&\bbE\left[\sup_{|t|\leq T}\Bigl|\|\psi_n\|_{H^1}^2 - \|\psi\|_{H^1}^2\Bigr|^{2p} \mathbbm 1_{\Omega_{n,A}}\right]\lesssim_{A,T,p} \epsilon_n^{\frac p4},
\end{align*}
whereas from \eqref{psi H1}, \eqref{uniform bounds for expectation in Xn1}, the Cauchy--Schwarz inequality, and \eqref{omega n A c}, we have
\begin{align*}
\bbE\left[\sup_{|t|\leq T}\Bigl|\|\psi_n\|_{H^1}^2 - \|\psi\|_{H^1}^2\Bigr|^{2p} \mathbbm 1_{\Omega_{n,A}^c}\right] &\lesssim_p\Biggl[ \Bigl[\bbE\|\psi_n\|_{C_TH^1}^{8p}\Bigr]^{\frac12} + \|\psi\|_{C_TH^1}^{4p}\Biggr]\bigl[\bbP(\Omega_{n,A}^c)\bigr]^{\frac12}\\
&\lesssim_{T,p} A^{-1} \bigl[1 + \|\psi_0\|_{H^1}^{4p+1}\bigr]\|\psi_0\|_{H^1}^{4p+1}.
\end{align*}
Combining these estimates, we obtain
\begin{align*}
\bbE\sup_{|t|\leq T}\Bigl|\|\psi_n\|_{H^1}^2 - \|\psi\|_{H^1}^2\Bigr|^{2p} &\lesssim_{T,p} C(A,T,p)\epsilon_n^{\frac p4} + A^{-1} \bigl[1 + \|\psi_0\|_{H^1}^{4p+1}\bigr]\|\psi_0\|_{H^1}^{4p+1},
\end{align*}
for some constant \(C(A,T,p)>0\). Taking \(n\to\infty\) and then $A\to \infty$ gives us
\eq{norm convergence}{
\lim_{n\to\infty}\bbE\sup_{|t|\leq T}\Bigl|\|\psi_n\|_{H^1}^2 - \|\psi\|_{H^1}^2\Bigr|^{2p} = 0.
}

We now estimate
\begin{align*}
\bbE \|\psi_n - \psi\|_{C_TH^1}^{4p} & \lesssim_p \bbE\sup_{|t|\leq T}\Bigl|\|\psi_n\|_{H^1}^2 - \|\psi\|_{H^1}^2\Bigr|^{2p} + \bbE\sup_{|t|\leq T}\bigl|\<\psi_n - \psi,\psi\>_{H^1}\bigr|^{2p},
\end{align*}
and hence \eqref{cvgce} holds whenever \(p\geq 4\). The case \(1\leq p<4\), then follows from applying H\"older's inequality to bound
\[
\bbE \|\psi_n - \psi\|_{C_TH^1}^p\leq \left[\bbE \|\psi_n - \psi\|_{C_TH^1}^{4p}\right]^{\frac 14},
\]
which completes the proof.
\epf

\section{Fluctuations}\label{s:fluctuations}
In this section we prove Theorem~\ref{t:cvgce in distrib}.
We start by considering the linear equation
\eq{linear toy}{
i\p_t u = - \p_x^2u + 4|\psi|^2u + 2\psi^2\bar u,
}
where \(\psi\in C(\R;H^1)\) is the solution of \eqref{NLS}. Our first lemma follows from classical arguments;  see, e.g., \cite{MR58861,MR54167}. As we have been unable to find the precise statement we require in the existing literature, for the reader's convenience we include a proof in Appendix~\ref{a:linear}:
\begin{lem}\label{l:toy}
Let \(\psi\in C(\R;H^1)\) be a solution of \eqref{NLS}. Then, for any \(|s|\leq 1\), any time \(\tau\in \R\), and any initial data \(u_0\in H^s\) there exists a unique solution \(u\in C(\R;H^s)\) of the Duhamel formulation,
\eq{linear prop}{
u(t) = e^{i(t-\tau)\p_x^2}u_0 - i\int_\tau^te^{i(t - \sigma)\p_x^2}\Bigl[ 4|\psi(\sigma)|^2u(\sigma) + 2\psi(\sigma)^2\overline{u(\sigma)}\Bigr]\,d\sigma,
}
of \eqref{linear toy} with \(u(\tau) = u_0\).

Moreover, if we define the real linear operator \(S(t,\tau)\colon H^s\to H^s\) by
\[
u(t) = S(t,\tau)u_0,
\]
the map \((t,\tau)\mapsto S(t,\tau)\) is continuous from \(\R^2\) to the space of bounded real linear operators on \(H^s\) endowed with the strong operator topology and satisfies the estimate
\eq{S op bound}{
\|S(t,\tau)\|_{H^s\to H^s}\leq \exp\Bigl(C|t - \tau|\bigl[1+\|\psi_0\|_{L^2}^2\bigr]\|\psi_0\|_{H^1}^2\Bigr),
}
for some \(C>0\).
\end{lem}

For the inhomogeneous equation,
\eq{inhomogeneous linear toy}{
i\p_t u = - \p_x^2u + 4|\psi|^2u + 2\psi^2\bar u + f,
}
we have following corollary that (for completeness) we also prove in Appendix~\ref{a:linear}:
\begin{cor}\label{c:inhomogeneous toy}
Let \(\psi\in C(\R;H^1)\) be a solution of \eqref{NLS}. Then, if \(|s|\leq 1\), \(\tau\in \R\), \(u_0\in H^s\), and \(f\in L_\loc^1(\R;H^s)\), the unique solution \(u\in C(\R;H^s)\) of the Duhamel formulation,
\eq{inhomogeneous Duhamel}{
u(t) = e^{i(t - \tau)\p_x^2}u_0 - i\int_\tau^te^{i(t - \sigma)\p_x^2}\Bigl[ 4|\psi(\sigma)|^2u(\sigma) + 2\psi(\sigma)^2\overline{u(\sigma)} + f(\sigma)\Bigr]\,d\sigma,
}
of \eqref{inhomogeneous linear toy} with initial data \(u(\tau) = u_0\) is given by
\eq{inhomogeneous form}{
u(t) = S(t,\tau)u_0 + \int_\tau^t S(t,\sigma)\bigl[\tfrac1if(\sigma)\bigr]\,d\sigma.
}
\end{cor}

Combining these results with Lemma~\ref{l:expect white noise} yields:
\begin{prop}\label{p:WP of lin-NLS}
If \(\psi\in C(\R;H^1)\) and \(s>\frac12\) then there almost surely exists a unique solution \(\varphi\in C(\R;H^{-s})\) of the Duhamel formulation, 
\eq{lin-NLS duhamel}{
\varphi(t) = - i\int_0^t e^{i(t - \tau)\p_x^2}\Bigl[4|\psi(\tau)|^2\varphi(\tau) + 2 \psi(\tau)^2\overline{\varphi(\tau)} + 2|\psi(\tau)|^2\psi(\tau)\,\xi\Bigr]\,d\tau,
}
of \eqref{lin-NLS} with initial data \(\varphi(0) = 0\).

Further, for each \(t\in \R\) the function \(\varphi(t)\) is a complex-valued Gaussian process with mean zero, covariance
\begin{align}
\bbE\bigl[\<f,\varphi\>\overline{\<g,\varphi\>}\bigr] &= \Bigl[\bigl\<\Re K_t^*f,\Re K_t^*g\bigr\> + \bigl<\Re K_t^*\bigl(\tfrac1i f\bigr),\Re K_t^*\bigl(\tfrac1i g\bigr)\bigr\>\Bigr]\label{covar}\\
&\quad + i\Bigl[\bigl\<\Re K_t^*\bigl(\tfrac1i f\bigr),\Re K_t^*g\bigr> - \bigl\<\Re K_t^*f, \Re K_t^*\bigl(\tfrac 1ig\bigr)\bigr\>\Bigr],\notag
\end{align}
and pseudo-covariance
\begin{align}
\bbE\bigl[\<f,\varphi\>\<g,\varphi\>\bigr] &= \Bigl[\bigl\<\Re K_t^*f,K_t^*g\bigr\> - \bigl<\Re K_t^*\bigl(\tfrac1i f\bigr),\Re K_t^*\bigl(\tfrac1i g\bigr)\bigr\>\Bigr]\label{pcovar}\\
&\quad + i\Bigl[\bigl\<\Re K_t^*\bigl(\tfrac1i f\bigr),\Re K_t^*g\bigr> + \bigl\<\Re K_t^*f,\Re K_t^*\bigl(\tfrac 1ig\bigr)\bigr\>\Bigr],\notag
\end{align}
for any complex-valued \(f,g\in H^s\), where \(K_t^*\colon H^s\to H^s\) is the dual of the real-linear map \(K_t\colon H^{-s}\to H^{-s}\) defined by
\eq{real IP}{
K_tf = \int_0^t S(t,\tau)\bigl[\tfrac{2}{i}|\psi(\tau)|^2\psi(\tau) f\bigr]\,d\tau
}
with respect to the pairing \(\Re\<f,g\>\) between \(f\in H^s\) and \(g\in H^{-s}\).
\end{prop}
\bpf Without loss of generality, we may assume that \(\frac12<s\leq 1\).

Given \(T>0\), we apply \eqref{product estimate}, the Cauchy--Schwarz inequality, and \eqref{expected white noise}, followed by the conservation of mass and \eqref{psi H1} to obtain
\begin{align*}
\bbE\bigl\||\psi|^2\psi\,\xi\bigr\|_{L_T^1H^{-s}}^2 &\lesssim T \|\psi\|_{C_TH^1}^4 \bbE \|\psi\xi\|_{L_T^2H^{-s}}^2\\
&\lesssim T\|\psi\|_{C_TH^1}^4\|\psi\|_{L_T^2L^2}^2\\
&\lesssim T^2\bigl[1 + \|\psi_0\|_{L^2}^4\bigr]\|\psi_0\|_{H^1}^4\|\psi_0\|_{L^2}^2.
\end{align*}

Choosing any sequence of times \(T_n\to\infty\), for all \(n\geq 1\) we have
\[
\bbE\bigl\||\psi|^2\psi\,\xi\bigr\|_{L_{T_n}^1H^{-s}}^2\lesssim T_n^2\bigl[1 + \|\psi_0\|_{L^2}^4\bigr]\|\psi_0\|_{H^1}^4\|\psi_0\|_{L^2}^2,
\]
and hence we almost surely have
\[
\bigl\||\psi|^2\psi\,\xi\bigr\|_{L_{T_n}^1H^{-s}}<\infty\qt{for all \(n\geq 1\)},
\]
so \(|\psi|^2\psi\xi\in L_\loc^1(\R;H^{-s})\) almost surely.

In the event that \(|\psi|^2\psi\,\xi\in L_\loc^1(\R;H^{-s})\), we now apply Corollary~\ref{c:inhomogeneous toy} to show that there exists a unique solution \(\varphi\in C(\R;H^{-s})\) of \eqref{lin-NLS} given by
\[
\varphi(t) = \int_0^t S(t,\tau)\Bigl[\tfrac{2}{i}|\psi(\tau)|^2\psi(\tau)\,\xi\Bigr]\,d\tau.
\]

Finally, fixing a time \(t\in \R\) and taking \(K_t\) to be defined as in \eqref{real IP}, we apply Lemma~\ref{l:toy}, \eqref{product estimate}, and \eqref{psi H1} to prove there exists a constant \(C>0\) so that
\[
\|K_t f\|_{H^{-s}}\leq\exp\Bigl(C|t|\bigl[1+\|\psi_0\|_{H^1}^2\bigr]\|\psi_0\|_{H^1}^2\Bigr)\|f\|_{H^{-s}}.
\]
By duality, we obtain
\[
\|K_t^*f\|_{H^s}\leq\exp\Bigl(C|t|\bigl[1+\|\psi_0\|_{H^1}^2\bigr]\|\psi_0\|_{H^1}^2\Bigr) \|f\|_{H^s},
\]
and hence
\begin{align*}
\bbE \exp\bigl(i\Re\<f,\varphi\>\bigr) &= \bbE\exp\bigl(i\<\Re K_t^*f,\xi\>\bigr)= \exp\bigl(-\tfrac12\|K_t^*f\|_{L^2}^2\bigr),
\end{align*}
where the second equality follows from \eqref{WN defn}. Using polarization identities, this suffices to prove that \(\varphi\) is a complex-valued, mean-zero Gaussian process with covariance \eqref{covar} and pseudo-covariance \eqref{pcovar}.
\epf

We now turn to the proof of Theorem~\ref{t:cvgce in distrib}. Our starting point is the following a priori estimate for solutions of \eqref{NLSn}:

\begin{lem}\label{l:ap psin again}
Let \(\frac12<s\leq1\), both \(\cN(0;Z_n)<\infty\) and \(\cN(0;Z_{n,s})<\infty\), and \(\psi_0\in Y_{n,s}^1\). Let \(\psi_n\in C(\R;X_n^1)\) be the solution of the Duhamel formulation \eqref{Duhamel NLSn} of \eqref{NLSn} with initial data \(\psi_n(0) = \psi_0\) obtained in Lemma~\ref{l:ap psin}. Then, \(\psi_n\in C(\R;Y_{n,s}^1)\) and for all \(T>0\) we have the estimate
\eq{refined ap for psin}{
\|\psi_n\|_{C_TY_{n,s}^1}^2\lesssim \<T\>\bigl[1 + \|\psi_0\|_{X_n^1}^2\bigr]\|\psi_0\|_{Y_{n,s}^1}^2.
}
\end{lem}
\bpf
As \(\psi_n\in C(\R;X_n^1)\) satisfies \eqref{Duhamel NLSn}, we may use \eqref{technical 2} to show that \(\psi_n\in C^1(\R;H^{-1})\) satisfies \eqref{NLSn} in \(H^{-1}\) for all \(t\in \R\). We may then argue precisely as in the proof of Lemma~\ref{l:ap psi}, replacing the estimate \eqref{psi H1} by \eqref{psin H1} to obtain \(\psi_n\in C(\R;L_{Z_{n,s}}^2)\) with the estimate
\[
\|\psi_n\|_{C_TL_{Z_{n,s}}^2}^2\lesssim \|\psi_0\|_{L_{Z_{n,s}}^2}^2 + T\bigl[1 + \|\psi_0\|_{X_n^1}^2\bigr]^{\frac12}\|\psi_0\|_{H^1}\|\psi_0\|_{L^2}.
\]
The estimate \eqref{refined ap for psin} now follows from combining this with \eqref{ap for psin}.
\epf

Recalling the definition \eqref{varphin}, we see that
\begin{align}
i\p_t\varphi_n &= - \p_x^2\varphi_n + 4|\psi|^2\varphi_n + 2\psi^2\bar\varphi_n + 2|\psi_n|^2\psi_n\tfrac{d\mu_n - 1}{\sqrt{\epsilon_n}}\label{varphin eq}\\
&\quad + 2|\psi_n - \psi|^2\varphi_n + 4\psi(\bar\psi_n - \bar \psi)\varphi_n + 2(\psi_n - \psi)\bar\psi\varphi_n,\notag
\end{align}
with initial data \(\varphi_n(0) = 0\).
By Corollary~\ref{c:inhomogeneous toy}, we may then write
\[
\varphi_n(t) = \cI_{1,n}(t) + \cI_{2,n}(t) + \cI_{3,n}(t),
\]
where
\begin{align}
\cI_{1,n}(t) &= \int_0^t S(t,\tau) \Bigl[\tfrac{2}{i}|\psi(\tau)|^2\psi(\tau)\tfrac{d\mu_n-1}{\sqrt{\epsilon_n}}\Bigr]\,d\tau,\label{I1n-def}\\
\cI_{2,n}(t) &= \int_0^t S(t,\tau) \Bigl[\tfrac{2}{i}\bigl[|\psi_n(\tau)|^2\psi_n(\tau) - |\psi(\tau)|^2\psi(\tau)\bigr]\tfrac{d\mu_n-1}{\sqrt{\epsilon_n}}\Bigr]\,d\tau,\label{I2n-def}\\
\cI_{3,n}(t) &= \int_0^t S(t,\tau)\Bigl[\tfrac{2}{i}\bigl|\psi_n(\tau) - \psi(\tau)\bigr|^2\varphi_n(\tau) + \tfrac{4}{i}\psi(\tau)\bigl[\bar\psi_n(\tau) - \bar \psi(\tau)\bigr]\varphi_n(\tau)\label{I3n-def}\\
&\qquad\qquad\qquad\qquad+ \tfrac{2}{i}\bigl[\psi_n(\tau) - \psi(\tau)\bigr]\bar\psi(\tau)\varphi_n(\tau)\Bigr]\,d\tau.\notag
\end{align}

In the spirit of our proof of Theorem~\ref{t:homogenization}, the main technical ingredient for the proof of Theorem~\ref{t:cvgce in distrib} is the following Proposition that collects several estimates under a boundedness assumption for \(\psi\) and \(\psi_n\); c.f.,~Proposition~\ref{p:the beast}.
\begin{prop}\label{p:technical nightmare}
Let \(\frac12<s\leq 1\). Suppose that \(0<\epsilon_n\leq 1\), \(\cN(0;Z_n)<\infty\), \(\cN(0;Z_{n,s})<\infty\), \(\psi_0\in Y_{n,s}^1\), and for some \(T,A>0\) we have
\eq{new A defn}{
\|\psi_n\|_{C_TY_{n,s}^1}^2 + \|\psi\|_{C_TY_{n,s}^1}^2 \leq A^2.
}
Then there exists a constant \(C = C(s)>0\) so that
\begin{align}
\|\varphi_n\|_{C_TH^{-s}} &\lesssim e^{CTA^2\opT}A,\label{varphin H-s}\\
\|\p_t\varphi_n\|_{C_TH^{-s-2}} &\lesssim e^{CTA^2\opT}(1+A^2)A,\label{varphin pt}\\
\|\cI_{2,n}\|_{C_TH^{-s}} + \|\cI_{3,n}\|_{C_TH^{-s}} &\lesssim e^{CTA^2\opT}A^{\frac12}\|\psi_n - \psi\|_{C_TH^1}^{\frac12}.\label{remainder terms}
\end{align}
where
\eq{triple norm}{
\opT = \sup_{|\tau|\leq |t|\leq T}\|S(t,\tau)\|_{H^{-s}\to H^{-s}}.
}

Further, for all \(R\geq 1\) we may bound
\begin{align}
&\|(1 - \chi_R)\varphi_n\|_{C_TH^{-s-1}}\label{varphin tight}\\
&\quad \lesssim e^{CTA^2\opT}AT\Bigl[A^{\frac32}\|(1 - \chi_R)\psi_n\|_{C_TH^1}^{\frac12} + A^{\frac32}\|(1 - \chi_R)\psi\|_{C_TH^1}^{\frac12} + \tfrac 1R\Bigr].\notag
\end{align}
\end{prop}
\bpf
Let us first note that Lemmas~\ref{l:ap psi} and~\ref{l:ap psin again} ensure that \(\varphi_n\in C(\R;Y_{n,s}^1)\) and hence \(\varphi_n\in C(\R;H^{-s})\).

Applying \eqref{technical 3} and the hypothesis \eqref{new A defn}, we may bound
\begin{align}
\|\cI_{1,n}\|_{C_TH^{-s}}&\lesssim T\opT\|\psi\|_{C_TY_{n,s}^1}^3 \lesssim T\opT A^3,\label{I1n}\\
\|\cI_{2,n}\|_{C_TH^{-s}}&\lesssim T\opT\Bigl[\|\psi\|_{C_TY_{n,s}^1}^{\frac52} + \|\psi_n\|_{C_TY_{n,s}^1}^{\frac52}\Bigr]\|\psi_n - \psi\|_{C_TH^1}^{\frac12}\notag\\
&\lesssim T\opT A^{\frac52}\|\psi_n - \psi\|_{C_TH^1}^{\frac12}.\label{I2n}
\end{align}
If we instead apply \eqref{product estimate} with \eqref{new A defn} then for \(|t|\leq T\) we have
\begin{align}
\|\cI_{3,n}(t)\|_{H^{-s}}&\lesssim \left|\int_0^t \opT \bigl[\|\psi_n\|_{C_TH^1} + \|\psi\|_{C_TH^1}\bigr]\|\psi_n - \psi\|_{C_TH^1}\|\varphi_n(\tau)\|_{H^{-s}}\,d\tau\right|\notag\\
&\lesssim \opT A \|\psi_n - \psi\|_{C_TH^1} \left|\int_0^t \|\varphi_n(\tau)\|_{H^{-s}}\,d\tau\right|\label{I3n}\\
&\lesssim T\opT A\|\psi_n - \psi\|_{C_TH^1}\|\varphi_n\|_{C_TH^{-s}}.\label{I3n'}
\end{align}

Combining the hypothesis \eqref{new A defn} with the estimates \eqref{I1n}, \eqref{I2n}, and \eqref{I3n}, for all \(|t|\leq T\) we may bound
\[
\|\varphi_n(t)\|_{H^{-s}}\lesssim \opT A^2\left|\int_0^t \|\varphi_n(\tau)\|_{H^{-s}}\,d\tau\right| + T \opT A^3.
\]
The estimate \eqref{varphin H-s} now follows from the integral form of Gronwall's inequality.

To prove \eqref{varphin pt}, we apply the estimates \eqref{product estimate} and \eqref{technical 3} to the expression \eqref{varphin eq} with the hypothesis \eqref{new A defn} and the estimate \eqref{varphin H-s}. The estimate \eqref{remainder terms} follows (after possibly increasing the size of \(C\)) from applying \eqref{new A defn} and \eqref{varphin H-s} to the estimates \eqref{I2n} and \eqref{I3n'}.

It remains to prove \eqref{varphin tight}. We first use that \(\varphi_n\) solves the Duhamel formulation of \eqref{varphin eq} with the identity
\[
[\chi_R,e^{it\p_x^2}] = -i\int_0^t e^{i(t-\tau)\p_x^2}\bigl(2\chi_R'\p_x + \chi_R''\bigr)e^{i\tau\p_x^2}\,d\tau,
\]
to write
\begin{align*}
(1 - \chi_R)\varphi_n(t) &= -2i\int_0^t e^{i(t-\tau)\p_x^2}\Bigl[(1 - \chi_R)\bigl(|\psi_n|^2\varphi_n + \psi_n\bar\varphi_n\psi + \varphi_n|\psi|^2\bigr)\Bigr]\,d\tau\\
&\quad - 2i\int_0^t e^{i(t-\tau)\p_x^2}\Bigl[(1 - \chi_R)|\psi_n|^2\psi_n\tfrac{d\mu_n-1}{\sqrt{\epsilon_n}}\Bigr]\,d\tau\\
&\quad + i\int_0^t e^{i(t-\tau)\p_x^2}\bigl[2\chi_R'\p_x\varphi_n(\tau) + \chi_R''\varphi_n(\tau)\bigr]\,d\tau.
\end{align*}
We may then apply \eqref{new A defn} and \eqref{varphin H-s} with the estimates \eqref{product estimate inf}, \eqref{product estimate}, and \eqref{technical 3} to bound
\begin{align*}
&\|(1 - \chi_R)\varphi_n\|_{C_TH^{-s-1}}\\
&\lesssim T \bigl[\|(1 - \chi_R)\psi_n\|_{C_TH^1} + \|(1 - \chi_R)\psi\|_{C_TH^1}\bigr]\bigl[\|\psi_n\|_{C_TH^1} + \|\psi\|_{C_TH^1}\bigr]\|\varphi_n\|_{C_TH^{-s}}\\
&\quad + T \|\psi_n\|_{C_TY_{n,s}^1}^2\|(1 - \chi_R)\psi_n\|_{C_TY_{n,s}^1}^{\frac12}\|(1 - \chi_R)\psi_n\|_{C_TH^1}^{\frac12} + \tfrac TR\|\varphi_n\|_{C_TH^{-s}}\\
&\lesssim \RHS{varphin tight},
\end{align*}
which completes the proof.
\epf

We now apply Prokhorov's Theorem with Lemma~\ref{l:compactness helper} and Proposition~\ref{p:technical nightmare} to obtain:

\begin{prop}\label{p:compactness}
For each \(T>0\) and \(s>\frac12\) the sequence of probability distributions \((\varphi_n)_*\bbP\) on \(C_TH^{-s}\) is precompact.
\end{prop}
\bpf
Without loss of generality, we assume that \(\frac12<s\leq1\).

We first prove that
\eq{R claim 1}{
\lim_{R\to\infty}\biggl[\tfrac1R + \sup_{n\geq 1}\bbE\|(1 - \chi_R)\psi_n\|_{C_TH^1}^{\frac12} + \|(1 - \chi_R)\psi\|_{C_TH^1}^{\frac12}\biggr] = 0.
}
For a fixed choice of \(n\geq 1\), Lemmas~\ref{l:ap psin} and \ref{l:ap psi} ensure that \(\psi_n - \psi\in C(\R;H^1)\) almost surely and hence we may apply Lemma~\ref{l:compactness helper} with \(\mc F = \{\psi_n - \psi\}\) to prove that
\[
\lim_{R\to\infty}\|(1 - \chi_R)(\psi_n - \psi)\|_{C_TH^1} = 0
\]
almost surely. From Theorem~\ref{t:homogenization}, we have \(\psi_n - \psi\in C_TH^1\) almost surely. In particular, for each \(n\geq 1\) we may use the Dominated Convergence Theorem with dominating function \(\|\psi_n - \psi\|_{C_TH^1}\) to obtain
\[
\lim_{R\to\infty}\bbE\|(1 - \chi_R)(\psi_n - \psi)\|_{C_TH^1} = 0.
\]
Next, we use \eqref{product estimate} to bound
\begin{align*}
&\lim_{R\to\infty}\sup_{n\geq 1}\bbE \|(1 - \chi_R)(\psi_n - \psi)\|_{C_TH^1}\\
&\qquad\lesssim \max_{n=1,\dots,N}\lim_{R\to\infty}\bbE \|(1 - \chi_R)(\psi_n - \psi)\|_{C_TH^1} + \sup_{n>N}\bbE \|\psi_n - \psi\|_{C_TH^1}\\
&\qquad\lesssim\sup_{n>N}\bbE \|\psi_n - \psi\|_{C_TH^1}.
\end{align*}
From \eqref{cvgce}, we may take \(N\to\infty\) to get
\[
\lim_{R\to\infty}\sup_{n\geq 1}\bbE\|(1 - \chi_R)(\psi_n - \psi)\|_{C_TH^1} = 0.
\]
Applying Jensen's inequality with \eqref{tight in H1} we then have
\begin{align*}
\LHS{R claim 1}&\lesssim \lim_{R\to\infty}\biggl[\sup_{n\geq 1}\bbE\|(1 - \chi_R)(\psi_n - \psi)\|_{C_TH^1} + \|(1 - \chi_R)\psi\|_{C_TH^1}\biggr]^{\frac12} = 0.
\end{align*}

We now set \(\sigma = \frac s2+\frac14\in(\frac12,s)\) and using \eqref{R claim 1} choose sequences of dyadic integers \(1\leq N_j,R_j\to\infty\) so that
\[
\sum_{j=1}^\infty 2^j \Bigl[N_j^{\sigma-s} + \tfrac1{R_j} + \sup_{n\geq 1}\bbE\|(1 - \chi_{R_j})\psi_n\|_{C_TH^1}^{\frac12} + \|(1 - \chi_{R_j})\psi\|_{C_TH^1}^{\frac12}\Bigr]<\infty.
\]
Given a constant \(M>0\), we then define the set \(\mc F_M\) to consist of all \(f\in C_TH^{-s}\) so that
\begin{align}
\|f\|_{C_TH^{-s}} + \|\p_tf\|_{C_TH^{-s-2}} &\leq M,\label{FM1}\\
\sum_{j=1}^\infty 2^j\|P_{>N_j}f\|_{C_TH^{-s}}&\leq M,\label{FM2}\\
\sum_{j=1}^\infty 2^j\|(1 - \chi_{R_j})f\|_{C_TH^{-s-1}} &\leq M,\label{FM3}
\end{align}
which is precompact in \(C_TH^{-s}\) by Corollary~\ref{c:FM}.

To complete the proof, we will prove that
\eq{prob measures are tight}{
\lim_{M\to\infty}\sup_{n\geq 1}(\varphi_n)_*\bbP(\mc F_M^c) = 0,
}
which shows that the sequence of probability measures \((\varphi_n)_*\bbP\) on \(C_TH^{-s}\) is tight. Prokhorov's Theorem \cite{MR84896,MR94838} now proves that this sequence of probability measures is precompact.

Given a constant \(A>0\), we take \(\Omega_{n,A}\) to be the event that \(\cN(0;Z_n)<\infty\), \(\cN(0;Z_{n,\sigma})<\infty\), \(\psi_0\in Y_{n,\sigma}^1\), and
\[
\|\psi_n\|_{C_TY_{n,\sigma}^1}^2 + \|\psi\|_{C_TY_{n,\sigma}^1}^2\leq A^2.
\]
As in the proof of Theorem~\ref{t:homogenization}, we may apply \eqref{ap for psi} and \eqref{refined ap for psin}, followed by \eqref{Xn1 exp} and \eqref{Yns1 exp} to prove that
\begin{align}
\sup_{n\geq 1}\bbE \Bigl[\|\psi_n\|_{C_TY_{n,\sigma}^1}^2 + \|\psi\|_{C_TY_{n,\sigma}^1}^2\Bigr] &\lesssim_{T,s}\bbE\Bigl[\bigl[1 + \|\psi_0\|_{X_n^1}^2\bigr]\|\psi_0\|_{Y_{n,\sigma}^1}^2\Bigr]\notag\\
&\lesssim_{T,s}\bigl[1 + \|\psi_0\|_{H^1}^2\bigr]\|\psi_0\|_{H^1}^2.\label{new uniform bounds for expectation}
\end{align}
Markov's inequality then yields the estimate
\begin{align}
\sup_{n\geq 1}\bbP\bigl[\Omega_{n,A}^c\bigr]&\leq A^{-2}\sup_{n\geq 1}\bbE \Bigl[\|\psi_n\|_{C_TY_{n,\sigma}^1}^2 + \|\psi\|_{C_TY_{n,\sigma}^1}^2\Bigr]\notag\\
&\lesssim_{T,s} A^{-2}\bigl[1 + \|\psi_0\|_{H^1}^2\bigr]\|\psi_0\|_{H^1}^2.\label{new omega n A c}
\end{align}

From \eqref{S op bound}, we have \(\opT\lesssim_{T,s,\psi_0}1\). In the event \(\Omega_{n,A}\), we may then apply the estimate \eqref{varphin H-s} to bound
\begin{align}
\bbE\Bigl[ \|P_{>N}\varphi_n\|_{C_TH^{-s}}\bbo_{\Omega_{n,A}}\Bigr] &\lesssim_s N^{\sigma-s}\bbE\Bigl[ \|P_{>N}\varphi_n\|_{C_TH^{-\sigma}}\bbo_{\Omega_{n,A}}\Bigr]\notag\\
&\lesssim_{A,T,s,\psi_0} N^{\sigma-s}.\label{compactness part a}
\end{align}
From the definition \eqref{Yns1}, we have \(\|f\|_{Y_{n,s}^1}\lesssim_s \|f\|_{Y_{n,\sigma}^1}\) and hence in the event \(\Omega_{n,A}\) we also have
\eq{new A defn cons}{
\|\psi_n\|_{C_TY_{n,s}^1}^2 + \|\psi\|_{C_TY_{n,s}^1}^2\lesssim_sA^2.
}
We may then apply \eqref{varphin H-s}, \eqref{varphin pt}, and \eqref{varphin tight} to obtain the estimates
\begin{align}
&\bbE\Bigl[\|\varphi_n\|_{C_TH^{-s}}\bbo_{\Omega_{n,A}}\Bigr] + \bbE\Bigl[\|\p_t\varphi_n\|_{C_TH^{-s-2}}\bbo_{\Omega_{n,A}}\Bigr] \lesssim_{A,T,s,\psi_0} 1,\label{compactness part c}\\
&\bbE\Bigl[\|(1 - \chi_R)\varphi_n\|_{C_TH^{-s-1}}\bbo_{\Omega_{n,A}}\Bigr]\notag\\
&\qquad\lesssim_{A,T,s,\psi_0} \tfrac1R + \bbE\|(1 - \chi_R)\psi_n\|_{C_TH^1}^{\frac12} + \|(1 - \chi_R)\psi\|_{C_TH^1}^{\frac12}.\label{compactness part b}
\end{align}
Finally, we apply \eqref{new omega n A c} and Markov's inequality with \eqref{compactness part a}, \eqref{compactness part c}, \eqref{compactness part b} and our choice of \(N_j,R_j\), to estimate
\begin{align*}
\sup_{n\geq 1}(\varphi_n)_*\bbP(K_M^c) &= \sup_{n\geq 1}\bbP\bigl[\varphi_n\not\in K_M\bigr]\\
&\leq\sup_{n\geq 1}\bbP\bigl[\Omega_{n,A}^c\bigr]\\
&\quad + \tfrac1M \sup_{n\geq 1}\biggl(\bbE\Bigl[\|\varphi_n\|_{C_TH^{-s}}\bbo_{\Omega_{n,A}}\Bigr] + \bbE\Bigl[\|\p_t\varphi_n\|_{C_TH^{-s-2}}\bbo_{\Omega_{n,A}}\Bigr]\biggr)\\
&\quad + \tfrac1M\sum_{j=1}^\infty 2^j \sup_{n\geq 1}\bbE\Bigl[\|P_{>N_j}\varphi_n\|_{C_TH^{-s}}\bbo_{\Omega_{n,A}}\Bigr]\\
&\quad + \tfrac1M\sum_{j=1}^\infty 2^j \sup_{n\geq 1}\bbE\Bigl[\|(1 - \chi_{R_j})\varphi_n\|_{C_TH^{-s}}\bbo_{\Omega_{n,A}}\Bigr]\\
&\lesssim_{T,s,\psi_0} \tfrac1{A^2} + \tfrac{C(A,T,s,\psi_0)}M,
\end{align*}
for some constant \(C(A,T,s,\psi_0)>0\). The identity \eqref{prob measures are tight} now follows from taking \(M\to\infty\) and then \(A\to\infty\).
\epf

We require one more lemma before we may prove Theorem~\ref{t:cvgce in distrib}:
\begin{lem}\label{l:key}
The maps
\begin{align*}
F&\mapsto \bbE\exp\biggl[\tfrac i{\sqrt{\epsilon_n}}\Re\biggl(\int_\R F\,d\mu_n -\int_\R F\,dx\biggr)\biggr],\\
F&\mapsto \bbE \exp\bigl[i\Re\<F,\xi\>\bigr],
\end{align*}
defined for \(F\in \Test\) have unique continuous extensions to \(F\in L^2\).

Moreover, if \(s\geq \frac16\) and \(F\in H^s\) then
\eq{CLT}{
\bbE\exp\biggl[\tfrac{i}{\sqrt{\epsilon_n}}\Re\biggl(\int_\R F\,d\mu_n -\int_\R F\,dx\biggr)\biggr] \to \bbE \exp\bigl[i\Re\<F,\xi\>\bigr]
}
as \(n\to\infty\).
\end{lem}
\bpf
Without loss of generality, we may assume that \(F\) is real-valued.

Let \(\Phi(z)\) be defined as in \eqref{new Phi expression} and \(x,y\in \R\).
By writing
\[
\Phi(ix) - \Phi(iy) = i(x-y)\int_0^1\Phi'\bigl(iy + ih(x-y)\bigr)\,dh,
\]
and using that
\[
|\Phi'(ix)| = \left|\int_0^\infty s(e^{-isx}-1)\,d\Lambda(s)\right|\lesssim |x|,
\]
we may bound
\[
\bigl|\Phi(ix) - \Phi(iy)\bigr| \lesssim \bigl(|x| + |y|\bigr)|x - y|.
\]
Similarly, recalling that \(\Phi(0) = \Phi'(0) = 0\) and \(\Phi''(0) = -1\), we may use Taylor's Theorem to bound
\[
\bigl|\Phi(ix) - \tfrac12x^2\bigr|\lesssim |x|^3,
\]
for all \(x\in \R\).

Given real-valued \(F\in L^2\), we then choose a real-valued sequence \(F_n\in \Test\) converging to \(F\) in \(L^2\). Recalling that \eqref{Laplace functional} applies to all \(f\in \Test(\R)\) satisfying \(\Re f>-\frac a{\epsilon_n}\) for \(a>0\) defined as in \eqref{Levy exponential bound}, we have
\begin{align*}
&\left|\bbE\exp\biggl[\tfrac i{\sqrt{\epsilon_n}}\biggl(\int_\R F_n\,d\mu_n -\int_\R F_n\,dx\biggr)\biggr] - \bbE\exp\biggl[\tfrac i{\sqrt{\epsilon_n}}\biggl(\int_\R F_m\,d\mu_n -\int_\R F_m\,dx\biggr)\biggr]\right|\\
&\quad = \left|\exp\left[-\tfrac1{\epsilon_n}\int_\R \Phi\bigl(-i\sqrt{\epsilon_n} F_n\bigr) - \Phi\bigl(-i\sqrt{\epsilon_n} F_m\bigr)\,dx\right] - 1\right|\\
&\quad\qquad\times \left|\exp\left[-\tfrac1{\epsilon_n}\int_\R \Phi\bigl(-i\sqrt{\epsilon_n} F_m\bigr)\,dx\right]\right|\\
&\quad\lesssim \exp\left[ C \bigl(\|F_n\|_{L^2}^2 + \|F_m\|_{L^2}^2\bigr)\right]\bigl(\|F_n\|_{L^2} + \|F_m\|_{L^2}\bigr)\|F_n - F_m\|_{L^2},
\end{align*}
for some \(C>0\).
Consequently, we may pass to the limit to see that (for real-valued \(F\)), the unique continuous extension of the map
\[
F\mapsto \bbE\exp\biggl[\tfrac i{\sqrt{\epsilon_n}}\biggl(\int_\R F\,d\mu_n -\int_\R F\,dx\biggr)\biggr]
\]
from \(\Test\) to \(L^2\) is given by
\[
\bbE\exp\biggl[\tfrac i{\sqrt{\epsilon_n}}\Re\biggl(\int_\R F\,d\mu_n -\int_\R F\,dx\biggr)\biggr] = \exp\left[-\tfrac1{\epsilon_n}\int_\R \Phi\bigl(-i\sqrt{\epsilon_n} F\bigr)\,dx\right].
\]

Similarly, for real-valued \(F\) the map
\[
F\mapsto \bbE \exp\bigl[i\<F,\xi\>\bigr] = \exp\bigl[-\tfrac12\|F\|_{L^2}^2\Bigr]
\]
has a unique continuous extension from \(\Test\) to \(L^2\).

Finally, if \(s\geq \frac13\) and \(F\in H^s\) is real-valued, we have
\begin{align*}
&\left|\bbE\exp\biggl[\tfrac i{\sqrt{\epsilon_n}}\biggl(\int_\R F\,d\mu_n -\int_\R F\,dx\biggr)\biggr] - \bbE \exp\bigl[i\<F,\xi\>\bigr]\right|\\
&\qquad = \left|\exp\left[-\tfrac1{\epsilon_n}\int_\R \Phi\bigl(-i\sqrt{\epsilon_n}F\bigr) - \tfrac{\epsilon_n}2|F|^2\,dx\right] - 1\right|\exp\left[-\tfrac12\|F\|_{L^2}^2\right]\\
&\qquad \lesssim \exp\left[C\sqrt{\epsilon_n} \|F\|_{L^3}^3\right]\sqrt{\epsilon_n}\|F\|_{L^3}^3.
\end{align*}
By Sobolev embedding we have \(\|F\|_{L^3}\lesssim \|F\|_{H^{\frac16}}\) and hence
\[
\lim_{n\to\infty}\bbE\exp\biggl[\tfrac i{\sqrt{\epsilon_n}}\biggl(\int_\R F\,d\mu_n -\int_\R F\,dx\biggr)\biggr] = \bbE \exp\bigl[i\<F,\xi\>\bigr],
\]
as claimed.
\epf

We are now finally in a position to complete the
\bpf[Proof of Theorem~\ref{t:cvgce in distrib}]
We recall (see, e.g., \cite{MR2814399}) that the characteristic functional of a random variable in  \(C_TH^{-s}\) uniquely determines its distribution.
In light of Proposition~\ref{p:compactness}, it remains to prove \eqref{fluctuation goal}.

Let \(\cT\colon C_TH^{-s}\to\C\) be a bounded linear operator and, recalling the definitions \eqref{I1n-def}--\eqref{I3n-def}, write
\begin{align}
\exp\bigl[i\Re\cT(\varphi_n)\bigr] &= \exp\bigl[i\Re\cT\bigl(\cI_{1,n}\bigr)\bigr]\Bigl[\exp\bigl[i\Re\cT\bigl(\cI_{2,n}+\cI_{3,n}\bigr)\bigr] - 1\Bigr]\label{exp i Re T}\\
&\quad + \exp\bigl[i\Re\cT\bigl(\cI_{1,n}\bigr)\bigr].\notag
\end{align}

For the first term on \RHS{exp i Re T}, we take \(\Omega_{n,A}\) to be the event that
\[
\|\psi_n\|_{C_TY_{n,s}^1}^2 + \|\psi\|_{C_TY_{n,s}^1}^2\leq A^2.
\]
As in the proof of Proposition~\ref{p:compactness}, we may first use \eqref{ap for psi}, \eqref{refined ap for psin}, \eqref{Xn1 exp}, and \eqref{Yns1 exp} to estimate
\[
\sup_{n\geq 1}\bbE \Bigl[\|\psi_n\|_{C_TY_{n,s}^1}^2 + \|\psi\|_{C_TY_{n,s}^1}^2\Bigr]\lesssim_{T,s,\psi_0}1,
\]
and then use Markov's inequality to obtain
\[
\sup_{n\geq 1}\bbP\bigl[\Omega_{n,A}^c\bigr]\lesssim_{T,s,\psi_0}A^{-2}.
\]
Combining this with \eqref{remainder terms} and Jensen's inequality, we may then bound
\begin{align*}
&\bbE\biggl|\exp\bigl[i\Re\cT\bigl(\cI_{1,n}\bigr)\bigr]\Bigl[\exp\bigl[i\Re\cT\bigl(\cI_{2,n}+\cI_{3,n}\bigr)\bigr] - 1\Bigr]\biggr|\\
&\qquad\lesssim \bbP\bigl[\Omega_{n,A}^c\bigr] + \bbE\biggl[ \|\cT\|_{(C_TH^{-s})^*}\Bigl[\|\cI_{2,n}\|_{C_TH^{-s}} + \|\cI_{3,n}\|_{C_TH^{-s}}\Bigr]\bbo_{\Omega_{n,A}}\biggr]\\
&\qquad\lesssim_{T,s,\psi_0}\tfrac1{A^2} + C(A,T,s,\psi_0,\cT)\Bigl[\bbE\|\psi_n - \psi\|_{C_TH^1}\Bigr]^{\frac12}.
\end{align*}
By first taking \(n\to\infty\) using \eqref{cvgce} and then \(A\to\infty\), we obtain
\eq{f g 1}{
\lim_{n\to\infty}\bbE\biggl|\exp\bigl[i\Re\cT\bigl(\cI_{1,n}\bigr)\bigr]\Bigl[\exp\bigl[i\Re\cT\bigl(\cI_{2,n}+\cI_{3,n}\bigr)\bigr] - 1\Bigr]\biggr| = 0.
}

For the remaining term on \RHS{exp i Re T}, we define \(q\in L_T^\infty H^s\) so that
\[
\Re\int_{-T}^T\int_\R\overline{q(\tau)}h(\tau)\,dx\,d\tau = \Re\cT\left(\int_0^tS(t,\tau)h(\tau)\,d\tau\right),
\]
for all \(h\in L_T^1H^{-s}\). We may then write
\[
\Re\cT(\cI_{1,n}) = \tfrac1{\sqrt{\epsilon_n}}\Re\left(\int_\R F\,d\mu_n - \int_\R F\,dx\right),
\]
where \eqref{product estimate} ensures that
\[
F = \tfrac2i\int_{-T}^T \overline{q(\tau)} |\psi(\tau)|^2\psi(\tau)\,d\tau\in H^s.
\]
Applying Lemma~\ref{l:key}, we obtain
\eq{f g 2}{
\lim_{n\to\infty}\bbE\exp\bigl[i\Re\cT\bigl(\cI_{1,n}\bigr)\bigr] = \bbE\exp\bigl[i\Re\<F,\xi\>\bigr] = \bbE\exp\bigl[i\Re\cT(\varphi)\bigr],
}
where \(\varphi\in C(\R;H^{-s})\) is the solution of \eqref{lin-NLS}.

The identity \eqref{fluctuation goal} now follows from combining \eqref{f g 1} and \eqref{f g 2}.
\epf

\begin{appendix}

\section{Solutions of \eqref{linear toy}}\label{a:linear}

In this appendix we prove Lemma~\ref{l:toy} and Corollary~\ref{c:inhomogeneous toy}.

\bpf[Proof of Lemma~\ref{l:toy}]
Given \(t\in \R\) and \(v\in H^s\), let
\[
L(t,v) = e^{-it\p_x^2}\Bigl[4|\psi(t)|^2 e^{it\p_x^2}v + 2\psi(t)^2 e^{-it\p_x^2}\overline{v}\Bigr].
\]
Employing the product estimate \eqref{product estimate} and recalling that the map \(t\mapsto e^{it\p_x^2}\) is continuous as a map from \(\R\) to the space of bounded linear operators on \(H^s\) endowed with the strong operator topology, we see that \(L\colon \R\times H^s\to H^s\) is continuous. Moreover, for any \(t\in \R\) and \(v,w\in H^s\) we may apply \eqref{product estimate} and \eqref{psi H1} to obtain the Lipschitz bound
\[
\|L(t,v) - L(t,w)\|_{H^s}\lesssim \bigl[1 + \|\psi_0\|_{L^2}^2\bigr]\|\psi_0\|_{H^1}^2 \|v - w\|_{H^s}.
\]
Standard ODE existence and uniqueness arguments then prove that for any \(u_0\in H^s\) there exists a unique global solution \(v\in C^1(\R;H^s)\) of
\eq{ODE}{
i\p_t v = L(t,v),
}
with initial data \(v(\tau) = e^{-i\tau\p_x^2}u_0\). Setting \(u = e^{it\p_x^2}v\in C(\R;H^s)\), we then see that there exists a unique global solution of \eqref{linear prop} satisfying \(u(\tau) = u_0\) with the estimate
\eq{ODE bound}{
\|u(t)\|_{H^s}\leq e^{C|t-\tau|\left[1 + \|\psi_0\|_{L^2}^2\right]\|\psi_0\|_{H^1}^2}\|u_0\|_{H^s},
}
for some constant \(C>0\).

The fact that \(S(t,\tau)\) is bounded, real linear, and satisfies \eqref{S op bound} follows immediately from the definition and \eqref{ODE bound}. 

To complete the proof, it remains to prove that for each \(u_0\in H^s\) the map \(\R^2\ni (t,\tau)\mapsto S(t,\tau)u_0\in H^s\) is continuous. To prove this, we will show that for any \(T>0\) we have
\begin{align}
\lim_{[-T,T]\ni \tau'\to \tau}\sup_{|t|\leq T}\|S(t,\tau')u_0 - S(t,\tau)u_0\|_{H^s} &= 0,\label{strong continuity 1}\\
\lim_{[-T,T]\ni t'\to t}\sup_{|\tau|\leq T}\|S(t',\tau)u_0 - S(t,\tau)u_0\|_{H^s} &= 0.\label{strong continuity 2}
\end{align}

Let us denote \(u_\tau(t) = S(t,\tau)u_0\). Using \eqref{linear prop}, we express
\begin{align*}
&u_{\tau'}(t) - u_\tau(t)\\
&\qquad= \bigl[e^{i(t-\tau')\p_x^2} - e^{i(t - \tau)\p_x^2}\bigr]u_0\\
&\qquad\quad - i\int_{\tau'}^t e^{i(t-\sigma)\p_x^2}\Bigl[4|\psi(\sigma)|^2\bigl[u_{\tau'}(\sigma) - u_\tau(\sigma)\bigr] + 2\psi(\sigma)^2\bigl[\overline{u_{\tau'}(\sigma)} - \overline{u_\tau(\sigma)}\bigr]\Bigr]\,d\sigma\\
&\qquad\quad - i\int_\tau^{\tau'} e^{i(t-\sigma)\p_x^2}\Bigl[4|\psi(\sigma)|^2u_\tau(\sigma) + 2\psi(\sigma)^2\overline{u_\tau(\sigma)}\Bigr]\,d\sigma.
\end{align*}
For \(t,\tau'\in [-T,T]\), we may then apply \eqref{product estimate} once more to estimate
\begin{align*}
\|u_{\tau'}(t) - u_\tau(t)\|_{H^s} &\lesssim \bigl\|\bigl[e^{-i\tau'\p_x^2} - e^{-i\tau\p_x^2}\bigr]u_0\bigr\|_{H^s}\\
&\quad + \|\psi\|_{C_TH^1}^2\left|\int_{\tau'}^t \|u_{\tau'}(\sigma) - u_\tau(\sigma)\|_{H^s}\,d\sigma\right|\\
&\quad + |\tau' - \tau|\|\psi\|_{C_TH^1}^2\|u_\tau\|_{C_TH^s}.
\end{align*}
Gronwall's inequality then yields the estimate
\begin{align*}
&\|u_{\tau'} - u_\tau\|_{C_TH^s}\\
&\qquad\lesssim e^{CT\|\psi\|_{C_TH^1}^2}\Bigl[\bigl\|\bigl[e^{-i\tau'\p_x^2} - e^{-i\tau\p_x^2}\bigr]u_0\bigr\|_{H^s} + |\tau' - \tau|\|\psi\|_{C_TH^1}^2\|u_\tau\|_{C_TH^s}\Bigr],
\end{align*}
for some constant \(C>0\), from which we obtain \eqref{strong continuity 1}.

For \eqref{strong continuity 2}, we again use \eqref{linear prop} to write
\begin{align*}
u_\tau(t') - u_\tau(t) &= \bigl[e^{i(t'-t)\p_x^2} - 1\bigr]u_\tau(t)\\
&\quad -i\int_t^{t'} e^{i(t'-\sigma)\p_x^2}\Bigl[4|\psi(\sigma)|^2 u_\tau(\sigma) + 2\psi(\sigma)^2\overline{u_\tau(\sigma)}\Bigr]\,d\sigma.
\end{align*}
Applying \eqref{product estimate} once again, we may bound
\[
\|u_\tau(t') - u_\tau(t)\|_{H^s}\lesssim \bigl\|\bigl[e^{i(t'-t)\p_x^2} - 1\bigr]u_\tau(t)\bigr\|_{H^s} + |t'-t| \|\psi\|_{C_TH^1}^2\|u_\tau\|_{C_TH^s}.
\]
From \eqref{strong continuity 1}, we see that the map \(\tau \mapsto u_\tau\) is continuous from \([-T,T]\to C_TH^s\). Consequently, the set
\(
\bigl\{u_\tau:|\tau|\leq T\bigr\}\subseteq C_TH^s
\)
is compact. Employing \eqref{compactness helper}, we then have
\[
\sup_{|\tau|\leq T}\|u_\tau\|_{C_TH^s}<\infty\qtq{and}\lim_{N\to\infty}\sup_{|\tau|\leq T}\|P_{>N}u_\tau\|_{C_TH^s} = 0.
\]
As a consequence, we may bound
\begin{align*}
\sup_{|\tau|\leq T}\|u_\tau(t') - u_\tau(t)\|_{H^s} &\lesssim |t' - t|N^2 \sup_{|\tau|\leq T}\|u_\tau\|_{C_TH^s} + \sup_{|\tau|\leq T}\|P_{>N}u_\tau\|_{C_TH^s}\\
&\quad + |t'-t| \|\psi\|_{C_TH^1}^2\sup_{|\tau|\leq T}\|u_\tau\|_{C_TH^s},
\end{align*}
and passing to the limit \(t'\to t\) and then \(N\to\infty\) gives us \eqref{strong continuity 2}.
\epf

\bpf[Proof of Corollary~\ref{c:inhomogeneous toy}]
Let us first note that uniqueness follows directly from the uniqueness of solutions to \eqref{linear toy}.

Turning to the problem of existence, let \(T>0\) and suppose that \(f\in C([\tau-T,\tau+T];H^s)\). Taking \(u\) to be defined as in \eqref{inhomogeneous form}, the strong continuity of the map \((t,\tau)\mapsto S(t,\tau)\) then ensures that \(u\in C([\tau-T,\tau+T];H^s)\).

From the proof of Lemma~\ref{l:toy}, for fixed \(\sigma\in[\tau-T,\tau+T]\) we have
\begin{align*}
i\p_t \Bigl[e^{-it\p_x^2}S(t,\tau)u_0\Bigr] &= e^{-it\p_x^2}\Bigl[4|\psi(t)|^2 S(t,\tau)u_0 + 2\psi(t)^2 \overline{S(t,\tau)u_0}\Bigr],\\
i\p_t \Bigl[e^{-it\p_x^2}S(t,\sigma)\bigl[\tfrac1if(\sigma)\bigr]\Bigr] &= e^{-it\p_x^2}\Bigl[4|\psi(t)|^2 S(t,\sigma)\bigl[\tfrac1if(\sigma)\bigr] + 2\psi(t)^2 \overline{S(t,\sigma)\bigl[\tfrac1if(\sigma)\bigr]}\Bigr],
\end{align*}
and hence
\[
i\p_t\Bigl[e^{-it\p_x^2}u(t)\Bigr] = e^{-it\p_x^2}\Bigl[4|\psi(t)|^2 u(t) + 2\psi(t)^2\overline{u(t)}\Bigr] + f(t).
\]
Integrating, we see that \(u\in C([\tau-T,\tau+T];H^s)\) solves \eqref{inhomogeneous Duhamel}.

The general case, where we only assume \(f\in L_\loc^1(\R;H^s)\) follows from fact that for all \(T>0\) the space \(C([\tau-T,\tau+T];H^s)\) is dense in \(L^1([\tau-T,\tau+T];H^s)\).
\epf

\section{Mollified measures}\label{app:mollification}
In this appendix we present analogs of Theorems~\ref{t:homogenization} and~\ref{t:cvgce in distrib} for the mollified model
\eq{NLSnh}{\tag{NLS$_n^h$}
i\p_t\psi_n^h = -\p_x^2\psi_n^h + 2|\psi_n^h|^2\psi_n^h\,\tfrac{d\mu_n^h}{dx},
}
where the density \(\frac{d\mu_n^h}{dx}\) is defined in \eqref{RND}. We state these as a single theorem:
\begin{thrm}\label{t:molly}
Let \(0<\epsilon_n\to0\) and \(\psi_0\in H^1\). For \(0<h\leq 1\), let \(\psi_n^h\in C(\R;H^1)\) be the almost surely defined solution of \eqref{NLSnh} with initial data \(\psi_n^h(0) = \psi_0\) and let \(\psi\in C(\R;H^1)\) be the solution of \eqref{NLS} with initial data \(\psi(0) = \psi_0\). Then:
\begin{enumerate}[(i)]
\item For all \(T>0\) and \(1\leq p<\infty\), we have
\eq{cvgce h}{
\lim_{n\to\infty}\sup_{0<h\leq 1}\bbE\sup_{|t|\leq T}\|\psi_n^h(t) - \psi(t)\|_{H^1}^p = 0.
}
\item For all \(T>0\) and \(s>\frac12\), if
\eq{varphin h}{
\varphi_n^h = \tfrac1{\sqrt{\epsilon_n}}(\psi_n^h - \psi),
}
and \(\varphi^h\in C(\R;H^{-s})\) is the almost surely defined solution of
\eq{lin-NLSh}{
i\p_t\varphi^h = - \p_x^2\varphi^h 
+ 4 |\psi|^2\varphi^h + 2\psi^2\bar\varphi^h + 2|\psi|^2\psi\,\xi^h
}
with initial data \(\varphi^h(0) = 0\), where \(\xi^h = \zeta^h*\xi\), then the sequence \(\varphi_n^h\) converges to \(\varphi^h\) in distribution as \(C\bigl([-T,T];H^{-s}\bigr)\)-valued random variables, uniformly for \(0<h\leq 1\).
\end{enumerate}
\end{thrm}

The proof of Theorem~\ref{t:molly} relies on a very simple observation, which already appeared in \cite[Lemma 2.2]{2024arXiv240501246H}. Let \(Z_n(k)\), \(Z_{n,s}(k)\) be defined as in \eqref{Xn1}, \eqref{Yns1}, respectively and let \(Z_n^h(k)\), \(Z_{n,s}^h(k)\) be the corresponding quantities with \(\mu_n\) replaced by \(\mu_n^h\). We may then bound
\eq{h variation}{
Z_n^h(k)\leq\sum_{\ell=k-1}^{k+1} Z_n(\ell)\qtq{and}Z_{n,s}^h(k)\lesssim_s \sum_{\ell=k-2}^{k+2} Z_{n,s}(\ell),
}
uniformly for \(0<h\leq 1\).
Now let \(X_n^{1,h}\), \(Y_{n,s}^{1,h}\) denote the spaces \(X_n^1\), \(Y_{n,s}^1\) with \(\mu_n\) replaced by \(\mu_n^h\). Employing \eqref{defn: cnk}, \eqref{sq sum} and \eqref{h variation}, it immediately follows that
\eq{uniform h}{
\sup_{0<h\leq 1}\|f\|_{X_n^{1,h}}\lesssim \|f\|_{X_n^1}\qtq{and}\sup_{0<h\leq 1}\|f\|_{Y_{n,s}^{1,h}}\lesssim \|f\|_{Y_{n,s}^1}.
}

With this in hand, let us now sketch the

\bpf[Proof of Theorem~\ref{t:molly}]
For part (i), we proceed precisely as in the proof of Theorem~\ref{t:homogenization} in Section~\ref{s:homogenization}, replacing \(\psi_n\) by \(\psi_n^h\), \(X_n^1\) by \(X_n^{1,h}\), and \(Y_{n,1}^1\) by \(Y_{n,1}^{1,h}\), with two exceptions. First, we use \eqref{uniform h} to replace \eqref{uniform bounds for expectation in Xn1} by
\begin{align}\label{new h-dependent uniform bounds}
\sup_{n\geq 1}\bbE \sup_{0<h\leq 1} \Bigl[\|\psi_n^h\|_{C_TX_n^{1,h}}^{2p} + \|\psi\|_{C_TY_{n,1}^{1,h}}^{2p}\Bigr]\
&\lesssim_{T,p}\bigl[1 + \|\psi_0\|_{H^1}^{2p}\bigr]\|\psi_0\|_{H^1}^{2p}.
\end{align}
Second, we now take \(\Omega_{n,A}\) to be the event that for all \(0<h\leq 1\) we have \(\cN(0;Z_n^h)<\infty\), \(\cN(0;Z_{n,1}^h)<\infty\), \(\psi_0\in Y_{n,1}^{1,h}\), and that \eqref{A defn} holds uniformly in \(h\), i.e.,
\[
\sup_{0<h\leq 1}\Bigl[\|\psi_n^h\|_{C_TX_n^{1,h}}^2 + \|\psi\|_{C_TY_{n,1}^{1,h}}\Bigr]\leq A^2.
\]
By applying \eqref{new h-dependent uniform bounds}, we see that the estimate \eqref{omega n A c} applies to our new version of \(\Omega_{n,A}\). This suffices to show that the remaining estimates hold uniformly in \(h\), giving us \eqref{cvgce h}.

The proof of part (ii) will follow a similar approach, although a few more modifications are required.

We first adapt Proposition~\ref{p:compactness} to prove that the set \(\bigl\{(\varphi_n^h)_*\bbP:n\geq 1,h\in(0,1]\bigr\}\) of probability distributions on \(C_TH^{-s}\) is precompact.

Rather than assuming \(\frac12<s\leq 1\) as in the proof of Proposition~\ref{p:compactness}, we will assume that \(\frac12<s<1\). From \cite[Lemma 3.3]{2024arXiv240501246H} we may bound
\[
\|(1-\chi_R)\psi_n^h\|_{C_TL^2}^2\lesssim \|(1-\chi_R)\psi_0\|_{L^2}^2 + \tfrac TR\|\psi_n^h\|_{C_TH^1}^2,
\]
uniformly for \(0<h\leq 1\) and \(n\geq 1\). We may then apply Bernstein's inequality for dyadic \(N\geq 1\) to bound
\begin{align*}
\|(1 - \chi_R)\psi_n^h\|_{C_TH^s} &\lesssim N^s \|(1-\chi_R)\psi_n^h\|_{C_TL^2} + N^{1-s}\|\psi_n^h\|_{C_TH^1}\\
&\lesssim N^s \|(1 - \chi_R)\psi_0\|_{L^2} + \Bigl[N^s\tfrac{\sqrt T}{\sqrt R} + N^{1-s}\Bigr]\|\psi_n^h\|_{C_TH^1},
\end{align*}
uniformly for \(n\geq 1\) and \(h\in(0,1]\). Consequently,
\begin{align*}
&\sup_{0<h\leq 1}\sup_{n\geq 1}\bbE\|(1 - \chi_R)\psi_n^h\|_{C_TH^s}\\
&\lesssim N^s \|(1 - \chi_R)\psi_0\|_{L^2} + \Bigl[N^s\tfrac{\sqrt T}{\sqrt R} + N^{1-s}\Bigr]\Bigl[\sup_{0<h\leq 1}\sup_{n\geq1}\bbE\|\psi_n^h-\psi\|_{C_TH^1} + \|\psi\|_{C_TH^1}\Bigr].
\end{align*}
From \eqref{cvgce h}, we have
\[
\sup_{0<h\leq 1}\sup_{n\geq1}\bbE\|\psi_n^h-\psi\|_{C_TH^1}<\infty,
\]
and hence we may take \(R\to\infty\) and then \(N\to\infty\) to obtain
\eq{R claim 1 h}{
\lim_{R\to\infty}\biggl[\tfrac1R + \sup_{0<h\leq 1}\sup_{n\geq 1}\bbE\|(1 - \chi_R)\psi_n^h\|_{C_TH^s}^{\frac1{2s}} + \|(1 - \chi_R)\psi\|_{C_TH^s}^{\frac1{2s}}\biggr] = 0.
}

We now take \(\mc F_M\) to be the precompact subset of \(C_TH^{-s}\) defined as in the proof of Proposition~\ref{p:compactness}, but with the dyadic sequences \(N_j,R_j\geq 1\) now chosen so that
\[
\sum_{j=1}^\infty 2^j \Bigl[N_j^{\sigma-s} + \tfrac1{R_j} + \sup_{0<h\leq 1}\sup_{n\geq 1}\bbE\|(1 - \chi_{R_j})\psi_n^h\|_{C_TH^s}^{\frac1{2s}} + \|(1 - \chi_{R_j})\psi\|_{C_TH^s}^{\frac1{2s}}\Bigr]<\infty,
\]
where \(\sigma = \frac s2+\frac14\in(\frac12,s)\) as before. As in part (i), we take \(\Omega_{n,A}\) to be the event that \(\cN(0;Z_n^h)<\infty\), \(\cN(0;Z_{n,\sigma}^h)<\infty\), and \(\psi_0\in Y_{n,\sigma}^{1,h}\) for all \(0<h\leq 1\), with the uniform estimate
\eq{another new A}{
\sup_{0<h\leq 1}\Bigl[\|\psi_n^h\|_{C_TY_{n,\sigma}^{1,h}}^2 + \|\psi\|_{C_TY_{n,\sigma}^{1,h}}^2\Bigr]\leq A^2.
}
Once again, we may use the estimate \eqref{uniform h} to prove that \eqref{new omega n A c} remains valid with this new definition of \(\Omega_{n,A}\).

The estimate \eqref{another new A} ensures that both \eqref{compactness part a} and \eqref{compactness part c}, with \(\varphi_n\) replaced by \(\varphi_n^h\), hold uniformly for \(0<h\leq 1\).

As \(\frac12<s<1\), we may use duality and the fact that \(H^s\) is an algebra to bound
\eq{new 23}{
\|fg\|_{H^{-s}}\lesssim \|f\|_{H^s}\|h\|_{H^{-s}}.
}
Combining this with the fractional Liebniz rule,
\[
\|fg\|_{H^s}\lesssim \|f\|_{L^\infty}\|g\|_{H^s} + \|f\|_{H^s}\|g\|_{L^\infty},
\]
the estimate \eqref{GN}, and the estimate
\[
\|f\|_{H^s}\lesssim \|f\|_{L^2}^{1-s}\|f\|_{H^1}^s,
\]
we may replace \eqref{technical 3} by
\begin{align}
\bigl\|f_1f_2g\tfrac1{\sqrt{\epsilon_n}}\bigl(\tfrac{d\mu_n^h}{dx} - 1\bigr)\bigr\|_{H^{-s}}&\lesssim \|f_1\|_{Y_{n,s}^{1,h}}\|f_2\|_{Y_{n,s}^1}\|g\|_{Y_{n,s}^{1,h}}^{1 - \frac1{2s}}\|g\|_{H^s}^{\frac1{2s}}.\label{technical 3000}
\end{align}
Replacing \eqref{product estimate} and \eqref{technical 3} by \eqref{new 23} and \eqref{technical 3000}, respectively, in the proof of \eqref{varphin tight} gives us
\begin{align*}
&\|(1 - \chi_R)\varphi_n^h\|_{C_TH^{-s-1}}\\
&\quad \lesssim e^{CTA^2\opT}AT\Bigl[A^{2-\frac1{2s}}\|(1 - \chi_R)\psi_n^h\|_{C_TH^s}^{\frac1{2s}} + A^{2-\frac1{2s}}\|(1 - \chi_R)\psi\|_{C_TH^s}^{\frac1{2s}} + \tfrac 1R\Bigr].
\end{align*}
Using this estimate instead of \eqref{varphin tight} in the proof of \eqref{compactness part b} yields the modification
\begin{align*}
&\bbE\Bigl[\|(1 - \chi_R)\varphi_n^h\|_{C_TH^{-s-1}}\bbo_{\Omega_{n,A}}\Bigr]\\
&\qquad\lesssim_{A,T,s,\psi_0} \tfrac1R + \bbE\|(1 - \chi_R)\psi_n^h\|_{C_TH^s}^{\frac1{2s}} + \|(1 - \chi_R)\psi\|_{C_TH^s}^{\frac1{2s}},
\end{align*}
which holds uniformly for \(0<h\leq 1\).

Combining these uniform estimates as in the proof of Proposition~\ref{p:compactness}, we arrive at the identity
\[
\lim_{M\to\infty}\sup_{0<h\leq 1}\sup_{n\geq 1}(\varphi_n^h)_*\bbP\bigl(\mc F_M^c\bigr) = 0.
\]
Prokhorov's Theorem now proves that the set \(\bigl\{(\varphi_n^h)_*\bbP:n\geq 1,h\in(0,1]\bigr\}\) of probability measures on \(C_TH^{-s}\) is precompact.

Next, we observe that the Laplace functional for the mollified measure \(\mu_n^h\) is
\[
\bbE\exp\left(-\int_\R f\,d\mu_n^h\right) = \bbE\left(-\int_\R \zeta^h*f\,d\mu_n\right),
\]
which suffices to prove that Lemma~\ref{l:key} applies with \(\mu_n\) replaced by \(\mu_n^h\) and \(\xi\) replaced by \(\xi^h\), where \eqref{CLT} holds uniformly for \(0<h\leq 1\).

Finally, a similar modification to the proof of Theorem~\ref{t:cvgce in distrib} in Section~\ref{s:fluctuations}, where we take \(\Omega_{n,A}\) to be the event that
\[
\sup_{0<h\leq 1}\Bigl[
\|\psi_n\|_{C_TY_{n,s}^{1,h}}^2 + \|\psi\|_{C_TY_{n,s}^{1,h}}^2\Bigr]\leq A^2,
\]
ensures that \eqref{fluctuation goal}, with \(\varphi_n\) replaced by \(\varphi_n^h\) and \(\varphi\) replaced by \(\varphi^h\), holds uniformly for \(0<h\leq1\).

To complete the proof, we proceed by contradiction. Suppose that there exists some \(\epsilon>0\), a bounded continuous function \(F\colon C([-T,T];H^{-s})\to \C\), and sequences \(n_j\to\infty\), \((h_j)\subseteq (0,1]\) so that
\eq{contradict me}{
\Bigl|\bbE F\bigl(\varphi_{n_j}^{h_j}\bigr) - \bbE F\bigl(\varphi^{h_j}\bigr)\Bigr|\geq \epsilon,
}
for all \(j\geq 1\).  Our modification of Proposition~\ref{p:compactness} ensures that, after passing to a subsequence, \(\varphi_{n_j}^{h_j}\) converges in distribution to a \(C_TH^{-s}\)-valued random variable \(\widetilde\varphi\).

After passing to a further subsequence, we may assume that \(h_j\to h\in[0,1]\). If we define the operator \(K^h = \<\p_x\>^{-s}\psi\,\zeta^h*\), we have
\[
\|K^{h_j} - K^h\|_{\I_2}^2 = \iint_{\R^2} \<\xi\>^{-2s}|\hat \psi(\xi-\eta)|^2 |\hat\zeta(h_j\eta)-\hat\zeta(h\eta)|^2\,d\xi\,d\eta,
\]
and an application of the Dominated Convergence Theorem proves that
\[
\lim_{j\to \infty}\|K^{h_j} - K^h\|_{\I_2} = 0.
\]
Arguing as in Lemma~\ref{l:expect white noise}, we then have
\[
\lim_{j\to\infty}\bbE\|\psi\xi^{h_j} - \psi\xi^h\|_{H^{-s}}^2 = 0,
\]
for any \(\psi\in L^2\). In particular, as in Proposition~\ref{p:WP of lin-NLS}, we may use Corollary~\ref{c:inhomogeneous toy} to bound
\[
\bbE\|\varphi^{h_j} - \varphi^h\|_{C_TH^{-s}}\lesssim \int_{-T}^T \opT \bbE \bigl\||\psi(\tau)|^2\psi(\tau) (\xi^{h_j} - \xi^h)\bigr\|_{H^{-s}}\,d\tau,
\]
where \(\opT\) is defined as in \eqref{triple norm}.
Another application of the Dominated Convergence Theorem then yields
\eq{varphi hj prob}{
\lim_{j\to\infty}\bbE\|\varphi^{h_j} - \varphi^h\|_{C_TH^{-s}} = 0.
}

Given a bounded linear functional \(\cT\colon C([-T,T];H^{-s})\to \C\), we now bound
\begin{align*}
&\biggl|\bbE e^{i\Re \cT(\widetilde\varphi)} - \bbE e^{i\Re \cT(\varphi^h)}\biggr|\\
&\qquad\leq\biggl|\bbE e^{i\Re \cT(\widetilde\varphi)} - \bbE e^{i\Re \cT(\varphi_{n_j}^{h_j})}\biggr| + \biggl|\bbE e^{i\Re \cT(\varphi_{n_j}^{h_j})} - \bbE e^{i\Re \cT(\varphi^{h_j})}\biggr|\\
&\quad \qquad + \biggl|\bbE e^{i\Re \cT(\varphi^{h_j})} - \bbE e^{i\Re \cT(\varphi^h)}\biggr|.
\end{align*}
Using that \(\varphi_{n_j}^{h_j}\) converges to \(\widetilde \varphi\) in distribution for the first summand,  applying our modification of \eqref{fluctuation goal} to the second summand, and using \eqref{varphi hj prob} for the third summand, we may pass to the limit as \(j\to\infty\) to show that \(\widetilde\varphi = \varphi^h\). Using \eqref{varphi hj prob} once again, we now see that for any bounded continuous function \(F\colon C([-T,T];H^{-s})\to \C\) we have
\[
\lim_{j\to\infty}\Bigl|\bbE F\bigl(\varphi_{n_j}^{h_j}\bigr) - \bbE F\bigl(\varphi^{h_j}\bigr)\Bigr| = 0,
\]
contradicting \eqref{contradict me}.
\epf

\end{appendix}

\bibliographystyle{habbrv}
\bibliography{refs}
\end{document}